\documentclass[10pt,a4paper,twoside]{amsart}

\usepackage{amsfonts, amssymb, amsmath, amsthm}
\usepackage{mathrsfs} 
\usepackage{latexsym}
\usepackage{enumerate}
\usepackage{verbatim}

\usepackage{url}
\usepackage{csquotes}
\DeclareUnicodeCharacter{0306}{\u i} 
\usepackage[%
    backend=bibtex,
    style=numeric-comp,
    citetracker=true,
    pagetracker=true,
    hyperref=true,
    backref=true,
    firstinits=true,
    bibencoding=utf8
    ]{biblatex}
\date{24 June 2025}

\usepackage[colorlinks=true]{hyperref}
\hypersetup{citecolor=blue, linkcolor=blue}

\newtheorem{thm}{Theorem}[section]
\newtheorem{lem}[thm]{Lemma}
\newtheorem{cor}[thm]{Corollary}

\newcommand{\proofsubsection}[1]{\par\smallskip\noindent{\sl #1.}}

\let\ve=\varepsilon

\def\gB{\mathfrak{B}}
\def\Bfrak{\mathfrak{B}}

\def\gP{\mathfrak{P}}

\def\Ocal{\mathcal{O}}

\def\1{1\!\!\!1}

\def\mode{\mathbin{\,\textrm{mod}^*}}
\def\mod{\mathbin{\,\textrm{mod}\,}}

\title[The trigonometric polynomial on sums of two squares]{The trigonometric polynomial on sums of two squares,
  an additive problem and generalisation}

\author{Olivier Ramar\'e}
\address[O. Ramar\'e]{CNRS/ Institut de Math\'ematiques de Marseille, Aix 
 Marseille Universit\'e, U.M.R. 7373, Site Sud, Campus de Luminy, Case 907, 
 13288 
 Marseille Cedex 9, France.}
\email{olivier.ramare@univ-amu.fr}

\author{G.K. Viswanadham}
\address[G. K. Viswanadham]{IISER Berhampur, Berhampur, Odisha, 760 010, India.}
\email{viswanadh@iiserbpr.ac.in}

\begin{document}

 \subjclass[2010]{Primary: 11L07, Secondary: 11M26, 11N25, 11P55}

 \keywords{Sums of two squares, Transference principle, Trigonometric polynomial}

\maketitle

\begin{abstract}
  Let $\Bfrak$ be the set of odd integers that are sums of two coprime
  squares. We prove that the trigonometric polynomial
  $S(\alpha;N)=\sum_{b\in\Bfrak,b\le N}e(b\alpha)$ satisfies
  \begin{equation*}
    \frac{S(\alpha; N)}{N/\sqrt{\log N}}\ll_{A,A'}
    \frac{1}{\varphi(q)}
    +
    \sqrt{\frac{q}{N}}(\log N)^{7}
    +\frac{1}{(\log N)^A}
  \end{equation*}
  for any $A,A'\ge0$ and when $(a,q)=1$ and $|q\alpha-a|\le (\log
  N)^{A'}/N$. We use this estimate together with a variant of the circle
  method influenced by Green and Tao's Transference Principle to
  obtain the number of representations of a large enough odd integer
  $N$ as a sum $b+b_1+b_2$, where $b\in\Bfrak$ while $b_1$
  (resp. $b_2$) belongs to a general subset $\Bfrak_1$
  (resp. $\Bfrak_2$) of $\Bfrak$ of relative positive density.
  We further show that the above bound is effective when $0\le A<1/2$.
\end{abstract}

{\small \tableofcontents}

\section{Introduction and results}

In 1922, the part III of the series Partitio Numerorum \cite{Hardy-Littlewood*22a} by {G.H.} {Hardy} \& {J.E.} {Littlewood} opened a new era saluted for instance by E.~Landau in his review for ZentralBlatt. In this paper, the authors developed the circle method introduced by G.H.~Hardy \& S.~Ramanujan in a series of papers, see~\cite{Hardy-Ramanujan*17-2}, to handle additive questions on primes. They proved in particular that, under the Generalized Riemann Hypothesis, every large enough odd integer is a sum of three primes, a result that was largely out of reach. Around the same time, V.~Brun in~\cite{Brun*19} initiated what was to become the combinatorial sieve. Influenced by both works, I.M.~Vinogradov simplified in~\cite{Vinogradov*37} Hardy \& Littlewood's method to use only the finite trigonometric polynomials on the primes, i.e. $\sum_{p\le N}e(p\alpha)$, and introduced a decomposition of this sum in what he called type~I or type~II sums (also referred to as linear and bilinear contribution), making possible to handle the minor arcs without appealing to the Generalized Riemann Hypothesis. This is a landmark in the study of primes. Later and simplifying considerably this initial work, P.X.~Gallagher 
in~\cite{Gallagher*68} and R.C.~Vaughan in~\cite{Vaughan*75} used convolution identities to obtain another type I/ type II decomposition. We finally mention that B.~Green \& T.~Tao in a series of papers around~\cite{Green-Tao*08-2} 
developed a dynamical understanding of this trigonometric polynomial; the restriction to primes still relies on combinatorial tools, as is surely easier to grasp from the paper \cite{Bourgain-Sarnak-Ziegler*12} by J.~Bourgain, P.~Sarnak \& T.~Ziegler (see \cite{Cafferata-Perelli-Zaccagnini*20} by M.~Cafferata, A.~Perelli \& A.~Zaccagnini).  
The survey paper \cite{Ramare*Vul*13-2} gives more details, but let us retain that two approaches compete to handle the trigonometric polynomials on the primes: via combinatorial tools or via convolution identities. Notice that the combinatorial tools are inherited from sieve techniques which rely heavily on non-negativity. 

A first aim of this paper is to try to disentangle the role of sieve from the one of convolution identities. The trigonometric polynomial~$S(\alpha;N)$ on sums of two squares is our trying ground. To stay in a clean sieve situation, we restrict our attention to the subset $\mathfrak{B}$ of odd integers that are sums of two coprime squares: these are simply the odd integers that are free of primes $\equiv3[4]$ and their characteristic function is written $n\mapsto b(n)$. The sieve dimension is just borderline for the sieve to be very efficient, see Chapter~14 of \cite{Friedlander-Iwaniec*10} by J.~Friedlander \& H.~Iwaniec. Convenient convolution identities are missing. 
We show that the initial Vinogradov idea (but not the one from Theorem~3, Chapter IX of \cite{Vinogradov*47}, see also \cite{Vinogradov*37d}) and the central identity of the combinatorial sieve may be used to produce adequate decompositions; we are able to cope with a more general case that we describe below, the two innovative steps being Theorem~\ref{Vino} and~\ref{AS}. We shall complete and improve on the latter in a forthcoming paper but decided here to stay focused on the trigonometric polynomial $S(\alpha;N)$;
we produce a finer bound for~$S(\alpha;N)$ when~$\alpha$ is close to a rational with a small denominator in Lemma~\ref{TrigoSmall}. We do not reach the maximal precision one gets for the primes, see \cite{Iwaniec-Ramare*09}, but miss it by a mere~$(\log\log N)^{\frac12+\ve}$. 
We end this study by a more classical asymptotic for $S(\alpha;N)$ when $\alpha$ is close to a rational with a small denominator. This work is supplemented by two applications, one to a ternary additive problem for which we prove an $L^p$-bound for $p>2$ for $S(\alpha;N)$ and another one to a quantitative measure of the equidistribution of $b\alpha$ when $\alpha$ is quadratic irrational and $b\in\gB$, for which we have to prove an average bound over $r\in\{1,2,\cdots,R\}$ for $S(r\alpha;N)$.
\bigskip

Before giving more details, we will introduce now a notation: $f=\Ocal^*(g)$ means that $|f|\le g$. We recall that $n\mapsto b(n)$ is the characteristic function of the set $\gB$ and 
we record that the function $\rho$ is defined in Eq.\,\eqref{defrho}.
\subsubsection*{Bilinear representations for characteristic functions
  of sifted sets}
To cover several cases as one, we introduce an
infinite set of primes $\gP$ that has a density $\kappa\in(0,1]$ in
the sense that
\begin{equation}
  \label{H1}
  \sum_{\substack{p\le X\\ p\in\gP}}\frac{\log p}{p}=\kappa\log X+\Ocal(1).
\end{equation}
When investigating the set~$\gB$ of odd integers that may be written
as sums of two coprime squares, $\gP$ is the set of primes congruent
to 3 modulo~4. Let us
notice that the elements of $\gB$ are all congruent to~1
modulo~4. This implies in particular that if a prime $p\equiv3[4]$
divides such an element, then at least another prime $\equiv 3[4]$ does so (which
may be the same, says $p$, meaning that $p^2$ divides our candidate).
As a conclusion, the odd integers below $N$ that are in~$\gB$ are the
integers below $N$ that are congruent to 1 modulo~4 and do not have
prime divisors $p\le \sqrt{N}$ congruent to~3 modulo~4.
 We can encompass also the case of
primitive Loeschian integers, i.e. integers~$n$ than can
be written as
$n=u^2+uv+v^2$ with $u$ and $v$ coprime by sieving the set of
integers
congruent to~1 modulo~3 by the set $\gP$
of primes congruent to~2 modulo~3.  As in
\cite{Fouvry-Levesque-Waldschmidt*18} by \'E.~Fouvry, C.~Levesque and
M.~Waldschmidt, we may also consider the set of integers that 
are both
odd primitive Gaussian integers and primitive Loeschian integers
in which case the
set $\gP$ is the set of primes not congruent to~1
mod~12. It is typographically expedient to define
\begin{equation}
  \label{defgPz}
  \gP(z)=\prod_{\substack{p<z\\ p\in\gP}}p
\end{equation}
and the special case
$P_{4;3}(z)=\prod_{\substack{p<z\\ p\equiv 3[4]}}p$.
Let us present our first identity that introduces few loss and, as may
be expected, leads to only a moderate saving. It results from a
combination of the initial method of I.~M.~Vinogradov together with
Rankin's trick. We state this initial result in a format close to the
one employed by H.~Iwaniec and E.~Kowalski in~\cite[Proposition
13.1]{Iwaniec-Kowalski*04}.
\begin{thm}
  \label{Vino}
  Let $2\le z\le D\le \sqrt{N}$ and $Z\ge z$. For any sequence of complex numbers
  $(u_n)$, we have
  \begin{equation*}
    \sum_{\substack{n\le N\\ (n,\gP(Z))=1}}u_n
    =
    \sum_{\substack{d|\gP(z)\\
        d\le D}}\mu(d)\sum_{n\equiv0[d]}u_n
    -\sum_{\substack{mp\le N\\  z\le p<
        Z}}
    \rho(m)\1_{p\in \gP}u_{mp}
    +\Ocal^*(E)
  \end{equation*}
  where
  \begin{equation}
  \label{defrho}
      \rho(m)=\frac{\1_{(m,\gP(z))=1}}
      {1+\#\{p'\in[z,Z]\cap\gP: p'|m\}}\in[0,1]
    \end{equation}
  and with
  \begin{equation}
    \label{defE}
    E = \sqrt{N\sum_{n}|u_n|^2}\biggl(
    \frac{1}{\sqrt{z}}
    +
    \exp\biggl(\frac{-\log D}{\log z}\biggr)(\log 3z)^{31}
    \biggr).
  \end{equation}
In case $|u_n|\leq 1$, the first term of the error term $E$ may be improved to $N/z$.
\end{thm}
This is proved in Section~\ref{ProofVino}.
Here is a consequence of it for the readers to visualize better.
\begin{cor}
  \label{Vinobn}
  Let $2\le z\le D\le \sqrt{N}$. For any sequence of complex numbers
  $(u_n)_{n\le N}$ supported on the integers congruent to~1 modulo~4, we have
  \begin{equation*}
    \sum_{n\le N}b(n)u_n
    =
    \sum_{\substack{d|P_{4;3}(z)\\
        d\le D}}\mu(d)\sum_{n\equiv0[d]}u_n
    -\sum_{\substack{mp\le N\\ z\le p\le
        \sqrt{N}}}
    \rho(m)\1_{p\equiv3[4]}u_{mp}
    +\Ocal^*(E)\; .
  \end{equation*}
   In case $|u_n|\leq 1$, the first term of the error term $E$ may be improved to $N/z$.
\end{cor}

  The non-negativity of $\rho(m)$ is noteworthy: it is an inheritance from the sieve background.
  The explicit form to the error is here to explicitate the
  dependencies in the parameters; 
  neither the constants, nor the power of logarithm have been optimized. 

Here is a more specific and interval oriented version of the previous result.
\begin{thm}
  \label{Vinobis}
  Let $10\le z\le D\le \sqrt{N}$ and $U\in[D, N]$ be another
  parameter. For any sequence of complex numbers
  $(u_n)$ supported in $[N-U,N]$ and such that
  $|u_n|\le 1$, Theorem~\ref{Vino} for $\gP(Z)=P_{4,3}(\sqrt{N})$ holds with $E$ replaced by
  $\Ocal(E')$ where
  \begin{equation*}
    E'=
    \frac{U}{z}
    +
    N\exp\biggl(\frac{-\log D}{\log z}\log\frac{\log D}{2\log z}\biggr)
    \sqrt{\log z}
    .
  \end{equation*}
\end{thm}
This is proved in Section~\ref{ProofVinobis}.
We complete our
weaponry with a somewhat less precise but much stronger bilinear
decomposition that combines the principle of Theorem~\ref{Vino} with
what we call the \emph{Simple sieve bilinear decomposition} contained in
Lemma~\ref{SieveBilinearSimple} and Theorem~\ref{PreciseAS}.
\begin{thm}
  \label{AS}
  For any set of parameters $\max(2,z)\le Z\le N$ and $2\le M_0\le M$,
  there exist sequences $(\alpha_\ell(t))$ and
  $(\beta_k(t))$ such that, for any sequence $(u_n)_{n\le N}$, we have
  \begin{multline*}
    \sum_{(n,\gP(Z))=1}\mkern-12mu u_n
    =\sum_{\substack{d<M\\d|\gP(z)}}\mu(d)\sum_{n\equiv 0[d]}u_n
    -\sum_{\substack{mp\le N\\ z\le p<
        Z}}
    \rho(m)\1_{p\in\gP}u_{mp}
    \\-
    \int_{0}^1
    \mkern-5mu\sum_{\substack{k\ell\ge M\\ M_0\le \ell\le M_0z}}\mkern-12mu
    \alpha_\ell(t)\beta_k(t)
    u_{k\ell}
    dt
    +
    \Ocal\biggl(
   \sqrt{ \frac{N\sum_{n}|u_n|^2}{z}}
    \biggr)
  \end{multline*}
  where
  $|\alpha_\ell(t)|\le 1$, $|\beta_k(t)|\le 16\tau_3(k)\log N$ and $\alpha_\ell(t)$ vanishes on indices~$\ell$ that do not divide
  $\gP(z)$, while $\rho(m)\in[0,1]$ is defined in Theorem~\ref{Vino}. In case~$|u_n|\le 1$, the error term may be improved to $\Ocal(N/z)$.
  \end{thm}
This is proved at the end of Section~\ref{UsingSimpleSieve}.
The function $\tau_3$ is classically the second Piltz function, namely $\tau_3(k)=\sum_{d_1d_2d_3=k}1$. We show in Theorem~\ref{thHS} that we may recover a version of the Harman sieve from the above decomposition. This version has a smaller admissible range for one of the variables.
\subsubsection*{Discrepancies between the case of primes and the case
  of sums of two squares}
Let us recall that the characteristic function of elements of $\gB$ (we
shall refer to them as \emph{odd primitive Gaussian integers})
is denoted by~$n\mapsto b(n)$.
One of our aim is to understand
the trigonometric polynomial
\begin{equation}
  \label{eq:2}
  S(\alpha; N)=\sum_{n\le N}b(n)e(n\alpha).
\end{equation}
Though one may think this polynomial over sums of
two squares is easier to handle than the one on the primes, the
sieving dimension being less, combinatorial identities \`a la
Gallagher-Vaughan are not available in this context.  Indeed we only
have
\begin{equation*}
  \sum_{\substack{d|n\\ d|P_{4;3}(\infty)}}\mu(d)
  =
  \begin{cases}
    1&\text{when $(n,P_{4;3}(\infty))=1$,}\\
    0&\text{otherwise,}
  \end{cases}
\end{equation*}
while, in the case of the primes, we may rely on
$\sum_{d|n}\mu(d)=\1_{n=1}$.

\subsubsection*{The trigonometric polynomial on sums of two squares}
The above machinery gives us access to sharp bounds for $S(\alpha;N)$ defined in~\eqref{eq:2}. By combining that with results on the distribution of sums of two squares in arithmetic progressions to small modulus, we reach the next theorem.
\begin{thm}
  \label{Trigo}
  Let $A\ge0$.
  Let $|\alpha-\frac{a}{q}|\le \frac{(\log N)^{2A+14}}{qN}$ with $(a,q)=1$ and $q\le N/(\log N)^{2A+14}$. We have
  \begin{equation*}
    \frac{S(\alpha; N)}{N/\sqrt{\log N}}\ll_A
    \frac{1}{\varphi(q)}
    +
    \sqrt{\frac{q}{N}}(\log N)^{7}
    +\frac{1}{(\log N)^A}
    .
  \end{equation*}
  The implied constant can be effectively computed when $A<1/2$.
\end{thm}
This is proved in Section~\ref{ProofTrigo}. See also Lemma~\ref{TrigoLargeII} and~\ref{TrigoSmall}.
There are two ranges of
estimation of $S(\alpha;N)$: when $q$ is small and when $q$ is large. 
We cover the large
range by using Theorem~\ref{AS} while the smaller range
of~$q$ is dealt with analytic methods, though
uneffectively. Effectivity can be recovered when $A<1/2$ by using
Corollary~\ref{Vinobn}.

\subsubsection*{Additive problems with subsequences of sums of two squares}
A theorem of Fermat in 1625 asserts that a positive integer $n$ can be
written as a sum of two coprime squares if and only if it is either odd and
not divisible by any prime $\equiv 3[4]$ or it is the double of such a
number. For the history of results on this subject, we refer the
readers to the Chapter~VI of Volume~II of the book~\cite{Dickson*71}
by L.~E.~Dickson.   In~\cite{Landau*08c}, E.~Landau proved
that
\begin{equation}
  \label{refC}
  B(N)=\sum_{n\le N}b(n)\sim
  \frac{CN}{\sqrt{\log N}},\quad
  C=\sqrt{2}\prod_{p\equiv3[4]}(1-1/p^2)^{-1/2}.
\end{equation}
Please note that such
integers are $\equiv 1[4]$. Our notation $b(n)$ is in accordance with
Section~14.3 of the reference book~\cite{Friedlander-Iwaniec*10} by
J.~Friedlander and H.~Iwaniec, see in particular Theorem 14.1 therein
for the value of $C$. Our aim is to study some additive questions
concerning subsets of $\Bfrak$. Their behaviour as essential
components and similar questions have already been treated in
\cite{Ramare-Ruzsa*01, Ramare*12-a} but only upper bounds were
required in such works. The half-dimensional sieve and sieve reversal principle may also
be put to work when it comes to the full sequence  $\Bfrak$, see
\cite[Corollary 2]{Iwaniec*76} by H.~Iwaniec, \cite{Greaves*76} by
G.~Greaves,  and
Section~14.5 of~\cite{Friedlander-Iwaniec*10}.
Queries requiring lower bounds or even
equality were successfully treated for rather general subsequences
of the primes in a series of works 
initiated by B.~Green in \cite{Green*05a}, followed by
\cite{Green-Tao*04-2, Green-Tao*10} by B.~Green and T.~Tao to cite but
a few of these results. We should mention specifically \cite{Shao*14a,
  Shao*14} by X.~Shao, \cite{Matomaki-Shao*17} by K.~Matom\"aki  and X.~Shao, \cite{Teravainen*18} by J.~Ter\"{a}v\"{a}inen or more recently \cite{Bauer*20} by C.~Bauer as
they are close to our main result.  In these pieces of work, the
\emph{transference principle}, as defined in \cite{Green-Tao*04} by
B.~Green and T.~Tao, is put to work. We are interested at avoiding
this method but to stick to the more precise circle method. We do so
below at the price of assuming that one variable is `more regular': 
X.~Shao \cite{Shao*14a} is able to handle sums of three primes, each belonging to a
dense subset of primes, whereas  we constrain one of our variables to belong to
the full set. As a result, we avoid the $W$-trick (but use restriction
estimates, see Lemma~\ref{L1} below) and get a close to asymptotic
expression for the required number of representations. 
The theorem we have in view is the following.
\begin{thm}
  \label{th2}
  Let $N\equiv3[4]$ be an integer, $K\in [2, \log\log N]$ be a
  parameter, and~$\mathcal{B}_1$ and~$\mathcal{B}_2$ be two sets of
  odd primitive Gaussian integers. With $M=P_{4,3}(K)$, we have
  \begin{equation*}
    \sum_{\substack{n+b_1+b_2=N\\ b_1\in\mathcal{B}_1, b_2\in\mathcal{B}_2}}\mkern-10mub(n)
    =
    \frac{CM/\varphi(M)}{\sqrt{\log N}}
    \mkern-20mu\sum_{\substack{b_1\in \mathcal{B}_1,
        b_2\in \mathcal{B}_2\\ b_1+b_2\le N\\
        (N-(b_1+b_2),M)=1.\\ }}
    \mkern-15mu
    1
    +\Ocal_{\mathcal{B}_1,\mathcal{B}_2}\biggl(
    \frac{N^2}{\sqrt{K}(\log N)^{3/2}}
    \biggr)
  \end{equation*}
  where $C$ is defined in~\eqref{refC}.
\end{thm}
This is proved at the end of Section~\ref{Proofth2}.
The variables $b_1$ and $b_2$ are controlled by the
restriction estimate recalled in Lemma~\ref{L1}. We may therefore 
replace $\gB_1$ and $\gB_2$ by any subsequence of relative positive density of a
\emph{sufficiently sifted set}, as defined in \cite{Ramare-Ruzsa*01},
e.g. the sequence of primes, or of primes~$p$ such that~$p+2$ is in
$\Bfrak$ as in \cite[Corollary 2]{Iwaniec*76} by H.~Iwaniec (see also
\cite{Indlekofer*74} by K.~H.~Indlekofer), or of Chen primes as in
\cite{Green-Tao*04}.
\subsubsection*{Sketching the proof of Theorem~\ref{th2}}
It is classical in the circle method to ``interpret" the singular series as the product of so-called ``local densities". We go one step further and provide
a geometrical meaning to this philosophy: we indeed approximate our function on the major arcs by a function that is periodical modulo~$M$, for some large modulus  $M$, that we call the \emph{local model}.

More precisely, the proof of Theorem~\ref{th2} goes by the circle method in which:
\begin{enumerate}
\item   The major arcs defined in \eqref{majorarcs} are arcs around $a/q$ for $a$ prime to~$q$ and $q\in P_{4,3}(K)\cup2P_{4,3}(K)\cup4P_{4,3}(K)$ for some parameter~$K$.
    \item We approximate one copy of $S(\alpha;N)$ by the trigonometric 
    polynomial $S^\flat_K(\alpha;N)$ of a local model of $\gB$ on the major arcs. 
    This approximation is the result of the comparison of
     Lemma~\ref{mainexpter} and~\ref{model}). 
\item The contribution of the minor arcs is shown to be negligible by
  combining two facts: (a) the $L^\ell$-inequality given in 
  Lemma~\ref{L1} for some $\ell>2$ and (b) the $L^\infty$-bound given 
  in Theorem~\ref{Trigo}. This idea is very much related with the work 
  of B.~Green \& T.~Tao, see Eq.\,(5.15) in~\cite{Green-Tao*04}. 
\item The expression with one copy of $S^\flat_K(\alpha;N)$ and 
  two copies of $S(\alpha;N)$ (though restricted to subsequences)
  is written back in terms of the coefficients, so $N=b_1+b_2+n$ with $n$ coming from the local model. The existence of such a writing amounts to a congruence condition of $N-b_1-b_2$ to some modulus~$M=4P_{4,3}(K)$, as in the first two lines on top of page~578 in~\cite{Ramare-Ruzsa*01}.
\end{enumerate}
The size of $M$, hence of $K$, is of utmost importance. The above scheme enables us to take~$K$ to be as small as a (large) constant if we rely on Theorem~\ref{Trigo}, i.e. 
if we allow analytical means. If we want to use only combinatorial means, 
so as to be able to generalize the proof, then Lemma~\ref{TrigoSmall}
 looses $\Theta=C\sqrt{\log\log N}\log\log\log N$, where $C$ is some positive
 constant; $K$ has to be larger than $\Theta^2$, implying that~$M$ is larger than any power of $\log N$. 
 
 Notice that the same proof holds if the three
variables are primes (or belonging to dense subsets of these for two of them) thanks, for the $L^\infty$-estimate, to
the Vinogradov estimates on the minor arcs coupled with `log'-free
estimates for the major arcs, as can be found in~\cite{Daboussi*96} by
H. Daboussi or in \cite{Iwaniec-Ramare*09} (see \cite{Ramare*13a}
for the latest best estimate) and, for the $L^\ell$-estimate, to the Bourgain-Green estimate found in Theorem 1.5 of \cite{Green*05}. 
On using the same setting as before, we may again take $K$ to be a constant.
This gives another Dirichlet
$L$-function-free proof of Vinogradov's three primes theorem, after
the one of X. Shao \cite{Shao*14}.  
\subsubsection*{A challenge}
Here is a result to measure the strength of our method in the middle range.
\begin{thm}
  \label{ValAtSqrtOf2}
  We have $S(\sqrt{2}; N)\ll N^{5/6}(\log N)^{13/2}$.
\end{thm}
This is proved in Section~\ref{ProofChallenge}.
Our process would give the same for primes, while the
Gallagher-Vaughan Identity immediately provides an exponent
$4/5$ rather than $5/6$. We could use any quadratic number rather than
$\sqrt{2}$ as the required property is that $\sqrt{2}$ has convergent
denominators of any given size, i.e. a substitute to
Lemma~\ref{ConVSqrt2}.

\subsubsection*{The trigonometric polynomial on the family of multiples}
The bilinear decomposition of Theorem~\ref{AS} enables us to handle sums of $e(b\alpha)$, but also of other oscillating
functions like $b\mapsto e(rb\alpha)$ for some integer~$r$.
Additional gain may be obtained by studying these functions on average
over a large range of values of $r$, and this leads to the next result. 
\begin{thm}
  \label{Family}
  Let $\alpha$ be a point of $\mathbb{R}/\mathbb{Z}$ given with a rational approximation $\alpha=\frac aq+\beta$ where $q^2|\beta|\le 1$. For $N\ge1$ and any $\ve>0$, we have
  \begin{equation*}
      \sum_{r\le R}|S(r\alpha;N)|\ll_\ve
      \frac{RN}{\sqrt{\log N}}
      \biggl(
      \frac{1}{\sqrt{q}}
      +\frac{\sqrt{q}}{\sqrt{RN}}
        +\frac{1}{R^{1/3}N^{1/6}}
        +\frac{R^{1/3}}{N^{1/3}}
      \biggr)
      N^\ve.
  \end{equation*}
\end{thm}
This is proved in Section~\ref{ProofFamily}. It can be compared 
with \cite[Theorem 3]{Vinogradov*54} by I.~M.~Vinogradov, 
or with \cite[Theorem 1]{Vaughan*77c} by R.~C.~Vaughan.
Specializing $R=1$ recovers the result of Lemma~\ref{TrigoLargeII}, save for the $(\log N)^7$ that is here degraded in $N^\ve$.

\subsubsection*{A Diophantine usage}
Here is a direct consequence of Theorem~\ref{Family}.
\begin{thm}
  \label{Irrat}
  Let $\alpha$ be a quadratic irrational and let $\lambda\in(0,1/4)$.
  The number of integers $b\in \gB$ with $b\le N$ and $\|b\alpha-\beta\|\le 1/N^\lambda$ is asymptotic
  to $2\gB(N)/N^\lambda$. 
\end{thm}
This is proved in Section~\ref{ProofIrrat}.
In the case of primes, R.~Vaughan in \cite{Vaughan*77c} obtained
a similar result. See also \cite[Theorem 2.2 and 3.2]{Harman*07} in the
book of G.~Harman.  More refined approaches by G.~Harman in
\cite{Harman*83, Harman*96} and by C.~H.~Jia in \cite{Jia*93b, Jia*00}
led to the exponent $\lambda=7/22$, though being unable to provide an
asymptotic formula in the case $\alpha$ quadratic irrational. In case
$\beta=0$, this exponent has been further reduced to $1/3$ after the
successive contributions of R.~Heath-Brown \& C.~H.~Jia in
\cite{HeathBrown-Jia*02} and of K.~Matom\"aki in \cite{Matomaki*09}.
Notice that, if we had only used Theorem~\ref{Trigo} rather than Theorem~\ref{Family}, the exponent $\lambda$ would be constrained to $\lambda<1/6$.

\subsubsection*{Thanks and acknowledgments}
Thanks are due to the referee for his/her very careful reading and for providing useful comments. The authors gratefully acknowledge funding from ReLaX, CNRS IRL 2024 and by NBHM,
and support and hospitality of IISER Berhampur and of I2M in Aix-Marseille University. 
\section{Proof of Theorem~\ref{Vino}}
\label{ProofVino}
Here is the first step.
\begin{lem}
\label{firstbase}
  Under the hypotheses of Theorem~\ref{Vino}, we have
  \begin{equation*}
  \sum_{\substack{n\le N\\ (n,\gP(Z))=1}}u_n
  =
  \sum_{\substack{n\le N\\ (n,\gP(z))=1}}u_n
  -\sum_{\substack{ z\le p\le Z\\ p\in\gP }}
  \sum_{\substack{m\le N/p}}
    \rho(m)u_{mp}
  +\Ocal^*\biggl(\sqrt{\frac{N}{z}\sum_{n}|u _n|^2}\biggr).
\end{equation*}
In case we have $|u_n|\le 1$, the error term may be improved to $\Ocal^*(N/z)$.
\end{lem}

\begin{proof}
  A moment of thought leads to our first formula:
\begin{equation*}
  \sum_{\substack{n\le N\\ (n,\gP(Z))=1}}u_n
  =
  \sum_{\substack{n\le N\\ (n,\gP(z))=1}}u_n
  -\sum_{\substack{ z\le p\le Z\\ p\in\gP }}
  \sum_{\substack{m\le N/p\\ (m,\gP(z))=1}}\frac{u_{pm}}{\omega^\sharp(pm)}
\end{equation*}
where
\begin{equation*}
  \omega^\sharp(pm)=\sum_{\substack{p'\in\gP\\ p'|pm\\ z\le p'\le Z}}1.
\end{equation*}
As this counting is done without multiplicities, we have
$\omega^\sharp(pm)=1+\omega(m)$ when $(m,p)=1$ and $\omega^\sharp(pm)=\omega^\sharp(m)$
otherwise. Hence
\begin{multline*}
   \sum_{\substack{n\le N\\ (n,\gP(Z))=1}}u_n
  =
  \sum_{\substack{n\le N\\ (n,\gP(z))=1}}u_n
  -\sum_{\substack{ z\le p\le Z\\ p\in\gP }}
  \sum_{\substack{m\le N/p\\ (m,\gP(z))=1}}
  \frac{u_{pm}}{1+\omega^\sharp(m)}
  \\-\sum_{\substack{ z\le p\le Z\\ p\in\gP}}
  \sum_{\substack{\ell\le N/p^2
      \\(\ell,\gP(z))=1}}\frac{u_{p^2\ell}}{(1+\omega^\sharp(p\ell))\omega^\sharp(p\ell)}.
\end{multline*}
We discard the last term in the following manner. First we notice that
\begin{equation}
  \label{step32}
  \sum_{\substack{ z\le p\le Z\\ p\in\gP}}
  \biggl|
  \sum_{\substack{\ell\le N/p^2
      \\(\ell,\gP(z))=1}}\frac{u_{p^2\ell}}{(1+\omega^\sharp(p\ell))\omega^\sharp(p\ell)}
  \biggr|
  \le \frac{1}2 
  \sum_{\substack{m\le N}}|u_{m}|\biggl(\sum_{\substack{ p^2|m\\ z\le p\le Z}}1\biggr).
\end{equation}
We apply the Cauchy-Schwarz inequality  to this last quantity and use
\begin{align*}
  \sum_{\substack{m\le N}}
  \biggl(\sum_{\substack{ p^2|m\\ z\le p\le Z}}1\biggr)^2
  &\le
  \sum_{z\le p_1 \neq p_2\le Z}\frac{N}{p_1^2p_2^2}
  +\sum_{z\le p_1\le Z}\frac{N}{p_1^2}
  \\&\le
  N\biggl(\frac{2}{z}\biggr)^2
  +\frac{2N}{z}\le \frac{4N}{z}.
\end{align*}
by using $\sum_{p\ge A}1/p^2\le 2/A$, valid for any $A\ge1$ and since we
may assume that $z\ge2$.
This proves our lemma with the first version of the error term. When
the sequence $(u_n)$ is bounded in absolute value by~1, we may proceed
more efficiently as in Theorem~26.2 of \cite{Moree-Ramare-Sedunova*22}
by P.~Moree, O.~Ramar\'e \& A.~Sedunova. Indeed we readily find that 
\begin{equation}
  \sum_{\substack{ z\le p\le Z\\ p\in\gP}}
  \biggl|
  \sum_{\substack{\ell\le N/p^2
      \\(\ell,\gP(z))=1}}\frac{u_{p^2\ell}}{(1+\omega^\sharp(p\ell))\omega^\sharp(p\ell)}
  \biggr|
  \le 
  \sum_{\substack{ z\le p\le Z\\ p\in\gP}}
  \frac{N}{2p^2}\le \frac{N}{z}.
\end{equation}
The proof is now to an end.
\end{proof}

\begin{proof}[Proof of Theorem~\ref{Vino}]
We first appeal to Lemma~\ref{firstbase}.
The first term may be studied by a straightforward usage of the
Erathostenes sieve by Legendre, i.e. we write
\begin{align*}
  \sum_{\substack{n\le N\\ (n,\gP(z))=1}}u_n
  &=\sum_{\substack{n\le N}}u_n\sum_{\substack{d|\gP(z)\\
      d|n}}\mu(d)
  \\&=
  \sum_{\substack{d|\gP(z)\\
  d\le D}}\mu(d)\sum_{n\equiv0[d]}u_n
  +
  \sum_{\substack{d|\gP(z)\\
  d> D}}\mu(d)\sum_{n\equiv0[d]}u_n.
\end{align*}
We treat the second summand in a similar manner to what we did to
bound the first error term, i.e. we write
\begin{equation}
  \label{step33}
  \sum_{\substack{d|\gP(z)\\
  d> D}}\mu^2(d)\sum_{n\equiv0[d]}|u_n|
  =
    \sum_{n\le N}|u_n|\sum_{\substack{d|\gP(z)\\
  d> D\\ d|n}}1
\end{equation}
to which we apply Cauchy's inequality to bound it above by the
squareroot of
\begin{equation*}
  \sum_{n\le N}|u_n|^2
  \sum_{n\le N}\biggl(\sum_{\substack{d|\gP(z)\\
      d> D\\ d|n}}1\biggr)^2.
\end{equation*}
We can handle the last sum, say $\Sigma$, via Rankin's trick:
\begin{align*}
  \Sigma
  &\le \sum_{\substack{D<d_1,d_2\le N, \\ [d_1,d_2]\le N\\
  d_1,d_2|\gP(z)}}
  \frac{N}{[d_1,d_2]}
  \le N\sum_{\substack{D<d\le N, \\
  d|\gP(z)}}
  \frac{3^{\omega(d)}}{d}
  \le N D^{-2\ve}\prod_{p\le z}\biggl(1+\frac{3p^{2\ve}}{p}\biggr)
  \\&\le N D^{-2\ve}\prod_{p\le z}\biggl(1+\frac{3p^{3\ve}}{p^{1+\ve}}\biggr)
  \le N D^{-2\ve} \zeta(1+\ve)^{3z^{3\ve}}
\end{align*}
for any positive $\ve$. We select $\ve=1/\log z$. Since
$\zeta(1+\ve)\le 1+1/\ve$, this finally amounts to
\begin{multline*}
  \sum_{\substack{n\le N\\ (n,\gP(Z))=1}}u_n
  =
  \sum_{\substack{d|\gP(z)\\
      d\le D}}\mu(d)\sum_{n\equiv0[d]}u_n
  -\sum_{\substack{mp\le N\\\ z\le p<
      Z\\ p\in\gP}}
  \rho(m)u_{mp}
  \\+\Ocal^*\biggl(
  \sqrt{N\sum_{n}|u_n|^2}\biggl(
  \frac{1}{\sqrt{z}}
  +
  \exp\biggl(\frac{-\log D}{\log z}\biggr)(1+\log z)^{31}
  \biggr)
  \biggr)\; .
\end{multline*}
This ends the proof of Theorem~\ref{Vino}.
\end{proof}

\section{Proof of Theorem~\ref{Vinobis}}
\label{ProofVinobis}
Before embarking into the proof, we recall a lemma due
to H.~Daboussi and J.~Rivat in this particular precise form. This is \cite[Lemma
4]{Daboussi-Rivat*01}, and is inspired by \cite[around page
81]{Elliott*79a} by P.D.T.A.~Elliott. 
\begin{lem}
\label{RankinDR}
Let $f$ be a non-negative multiplicative function, and $z\ge2$ be some
real parameter. We define
\begin{equation*}
    S=\sum_{p<z}\frac{f(p)}{1+f(p)}\log p\;.
\end{equation*}
We assume $S>0$ and write $K(t)=\log t-1+1/t$ for $t\ge1$. Then for
any $y$ such that $\log y\ge S$, we have
\begin{equation*}
  \sum_{\substack{d\ge y\\ d|P(z)}}\mu^2(d)f(d)
  \le \prod_{p<z}(1+f(p))\exp\biggl(-\frac{\log y}{\log z}
  K\biggl(\frac{\log y}{S}\biggr)\biggr)\;.
\end{equation*}
When $\log y\ge 7S$, we have $K\bigl((\log y)/{S}\bigr)\ge1$.
The function $t\mapsto K(t)$ is non-decreasing when
$t\ge1$.
\end{lem}

Here is the consequence of Lemma~\ref{RankinDR} that we are interested in.
\begin{lem}
  \label{RTBrun}
 When $z\ge 10$ and $2\log D\ge \log z$, we have 
  \begin{equation*}
    \sum_{D<d|P_{4;3}(z)}\frac{\mu^2(d)}{d}
    \ll \sqrt{\log z}
    \exp\biggl(\frac{-\log D}{\log z}\log\frac{\log D}{2\log z}\biggr)\; .
  \end{equation*}
\end{lem}

\begin{proof}
  We apply Lemma~\ref{RankinDR} with the multiplicative function $f$
  defined on primes by
  \begin{equation}
    \label{eq:6}
    f(p)=
    \begin{cases}
     1/p &\text{when $p\equiv3[4]$}\\
      0&\text{otherwise.}
    \end{cases}
  \end{equation}
  The values on prime powers are irrelevant, so we set $f(p^k)=0$ when
  $k\ge2$ to be semantically correct. We use Lemma~\ref{mertensp34} to
  compute that
  \begin{equation*}
    S=\sum_{p<z}\frac{f(p)\log p}{1+f(p)}
    =\sum_{\substack{p<z\\ p\equiv 3[4]}}\frac{\log p}{p+1}
    \le \tfrac12\log z -\tfrac14\log 3
    \le \frac{\log z}{2} .
  \end{equation*}
  
  By monotonicity, we get, when $z\ge 10$ and $2\log D\ge \log z$
  \begin{align*}
    K\biggl(\frac{\log D}{S}\biggr)
    &\ge K\biggl(\frac{2\log D}{\log z}\biggr)
    \ge \log\frac{2\log D}{\log z}-1
      \\&\ge  \log\frac{\log D}{2\log z}.
  \end{align*}
  This concludes the lemma.
\end{proof}

\begin{proof}[Proof of Theorem~\ref{Vinobis}]
  We follow the proof of Theorem~\ref{Vino}. Since $|u_n|\leq 1$, we get that the first error term is bounded by $\frac{U}{z}$.   The next
  discrepancies appear after~\eqref{step33}. We again use the fact
  that $|u_n|\le 1$ to bound above $|\sum_{n:d|n}u_n|$ by $N/d$. We
  then appeal to Lemma~\ref{RTBrun} to get the required result.
\end{proof}

\section{A Simple Sieve Identity}
The main result of this section is Lemma~\ref{SieveBilinearSimple} to which we provide two proofs. Let us start by a more general version of it. 

\begin{lem}
  \label{SieveBilinear}
  Let~$\nabla$ be a function such that (1) $\nabla$ takes only the values~0 and~1, (2) the set on which $\nabla(d)=1$ is divisor-closed and (3) $\nabla(1)=1$.
 Then for any sequence~$\psi$, we have
  \begin{align*}
    \sum_{d|\gP(z)}\mu(d)\psi(d)
    =&\sum_{d|\gP(z)}\mu(d)\psi(d)\nabla(d)\notag
       +\sum_{\ell|\gP(z)}
       \mu(\ell)\psi(\ell)\overline{\nabla}(\ell)
    \\&
    -\sum_{\ell|\gP(z)}\overline{\nabla}(\ell)\mu(\ell)
    \sum_{\substack{p'<P^{-}(\ell)\\ p'\in\gP}}
    \sum_{\substack{d|\gP( p')}}
    \mu(d)\psi(\ell p'd)
  \end{align*}
  where
  \begin{equation}
  \label{eq:12}
  \overline{\nabla}(\ell)
  =
  \begin{cases}
    0&\text{when $\ell=1$}\\
    \nabla(\ell/P^{-}(\ell))-\nabla(\ell)
    &\text{when $\ell>1$}.
  \end{cases}
\end{equation}
\end{lem}
We recall that $P^{-}(m)$
denotes the least prime factor of the integer~$m$.
We shall see that the binding condition $p'< P^{-}(\ell)$ is easy to remove.
On specializing $\nabla$ (and some shuffling), we get the next simpler identity.
\begin{lem}
  \label{SieveBilinearSimple}
  For any sequence $\psi$ and $M\ge M_0\ge 1$, we have
  \begin{multline*}
    \sum_{d|\gP(z)}\mu(d)\psi(d)
    =
    \sum_{\substack{d|\gP(z)\\ d< M}}\mu(d)\psi(d)
    +
    \sum_{\substack{\ell|\gP(z)\\ M\le \ell<M_0 P^{-}(\ell)}}
    \mu(\ell)\psi(\ell)
    \\-\sum_{\substack{\ell|\gP(z)\\ M_0\le \ell<M_0 P^{-}(\ell)}}
    \mu(\ell)
    \sum_{\substack{p'< P^{-}(\ell)\\ p'\in\gP(z)}}
    \sum_{\substack{d|\gP(p')\\ \ell dp'\ge M}}\mu(d) 
    \psi(\ell dp').
  \end{multline*}
\end{lem}
When specializing $M=M_0$, this identity takes a simpler look, which
is directly inferred from Lemma~\ref{SieveBilinear}. The distinction
between the parameters $M$ and $M_0$ may be useful.

\begin{proof}[Proof of Lemma~\ref{SieveBilinear}]
  Let us recall the Fundamental Sieve Identity.  This is, in raw form,
  Eq.~(1.8) of \cite[Chapter 2]{Halberstam-Richert*74} by
  H.~Halberstam and H.-E.~Richert, and, more explicitly, Lemma~2.1 of
  \cite{Diamond-Halberstam-Richert*88} of H.~Diamond, H.~Halberstam
  and H.-E.~Richert (reproduced with some notational changes in
  Lemma~2 of \cite{Ford-Halberstam*00} by K.~Ford and
  H.~Halberstam). It reads:
  \begin{equation}
    \label{FSI}
    \sum_{d|\gP(z)}\mu(d)\psi(d)
    =\sum_{d|\gP(z)}\mu(d)\psi(d)\nabla(d)\notag
    +\sum_{\ell|\gP(z)}\overline{\nabla}(\ell)
    \sum_{\substack{\delta|\gP( P^-(\ell))}}
    \mu(\ell\delta)\psi(\ell\delta).
  \end{equation}
  We have used Lemma~2 of
  \cite{Ford-Halberstam*00} with $h(d)=\mu(d)\psi(d)$ and $D=\gP(z)$; then
   in the second sum, $d$ and $t$ are replaced by~$\ell$ and $\delta$ respectively.
  Notice that $\mu(\ell\delta)=\mu(\ell)\mu(\delta)$.
  As $\delta$ and $\ell$ are linked in a tricky manner, we use the Buchstab
  Identity (see Eq. (1.9) of \cite[Chapter~2]{Halberstam-Richert*74} or~\eqref{BI} below) to write
  \begin{align*}
    \sum_{d|\gP(z)}\mu(d)\psi(d)
    =&\sum_{d|\gP(z)}\mu(d)\psi(d)\nabla(d)\notag
       +\sum_{\ell|\gP(z)}\overline{\nabla}(\ell)
       \mu(\ell)\psi(\ell)
    \\&
    -\sum_{\ell|\gP(z)}\overline{\nabla}(\ell)\mu(\ell)
    \sum_{\substack{p'<P^{-}(\ell)\\ p'\in\gP}}
    \sum_{\substack{d|\gP( p')}}
    \mu(d)\psi(\ell p'd)
  \end{align*}
  which is exactly what we announced.
\end{proof}
\begin{proof}[Proof of Lemma~\ref{SieveBilinearSimple}]
  When $M=M_0$, the selection~$\nabla(d)=\1_{d<M}$ in Lemma~\ref{SieveBilinear}
  yields Lemma~\ref{SieveBilinearSimple}. For the general case, use
  Lemma~\ref{SieveBilinear}
  on
  $\psi(d)\1_{d\ge M}$ and $\nabla(d)=\1_{d<M_0}$. The reader will swiftly
  recover the claimed identity.
\end{proof}

Lemma~\ref{SieveBilinearSimple} is at the heart of our proof, and it is
methodologically interesting to note that it may be derived directly from the next decomposition. 
\begin{lem}
  \label{HarmanDec}
  Let $M_0\ge 1$ be a real parameter and let $d\ge M_0$ be a squarefree
  integer. Then there exists a unique couple $(\delta,\ell)$ such that
  $d=\delta\ell$, $\delta|\gP(P^-(\ell))$ and
  $M_0\le \ell< M_0P^{-}(\ell)$.
\end{lem}
\begin{proof}
  Let us write $d=p_1p_2\cdots p_r$ where $p_1>p_2>\cdots>p_r$.
  Since $M_0\leq d$, there exists a $t_0$ such that
  $p_1\cdots p_{t_0-1}<M_0\leq p_1\cdots p_{t_0}$. We choose
  $\ell=p_1\cdots p_{t_0-1}p_{t_0}$ and $\delta=p_{t_0+1}\cdots
  p_r$. Clearly $M_0\leq \ell=p_1\cdots p_{t_0}$ and
  $M_0p_{t_0}> p_1\ldots p_{t_0}=\ell$. Suppose there is another
  decomposition of $\ell$ as $\ell=\delta_1\ell_1$ with the same properties.
  Then there exists an $s\leq r$ such that $\ell_1=p_1\cdots p_s$ and
  $\delta_1=p_{s+1}\cdots p_r$. We have
  \[
  p_1p_2\ldots p_{t_0-1}<M_0\leq p_1\cdots p_s
  \] 
  which implies $t_0\leq s$. On the other hand
  $\ell_1=p_1\cdots p_s<M_0p_s<p_1\cdots p_{t_0}p_s$. If
  $s\geq t_0+1$, we get that $\prod_{i=t_0+1}^{s-1}p_i<1$, a
  contradiction.
\end{proof}

\begin{proof}[Second proof of Lemma~\ref{SieveBilinearSimple}]
  By Lemma~\ref{HarmanDec}, we find that
  \begin{equation}
  \label{annoted}
    \sum_{d|\gP(z)}\mu(d)\psi(d)
    =
    \sum_{\substack{d|\gP(z)\\ d< M}}\mu(d)\psi(d)
    +
    \sum_{\substack{\ell|\gP(z)\\ M_0\le \ell<M_0 P^{-}(\ell)}}
    \mu(\ell)
    \sum_{\substack{\delta|\gP(P^{-}(\ell))\\ \delta\ell \ge M}}\mu(\delta) 
    \psi(\ell \delta).
  \end{equation}
  Once this identity is established, the binding condition  $\delta|\gP(P^{-}(\ell))$ on  $\delta$ and $\ell$ can be handled using the Buchstab Identity:
  \begin{equation}
    \label{BI}
    \sum_{\delta|\gP(P^{-}(\ell))}\mu(\delta) 
    \1_{\ell\delta\ge M}\psi(\ell \delta)
    =
    \1_{\ell\ge M}\psi(\ell)
    -\sum_{\substack{p'< P^{-}(\ell)\\ p'\in\gP}}
    \sum_{d|\gP(p')}\mu(d)
    \1_{\ell dp'\geq M}\psi(\ell dp').
  \end{equation}
  This ends the proof.
\end{proof}

\section{Using the Simple Sieve Identity}
\label{UsingSimpleSieve}
We use Fourier analysis to separate $\ell$ and
$p'$ in Lemma~\ref{SieveBilinear}.
\begin{lem}
  \label{VFA}
  Let $u,v>0$ with $v\neq u$, and let $T\ge1$. We have
  \begin{equation*}
    \frac{4}{\pi}\int_{1/T}^T\sin^2(vt/2)\frac{\sin ut}{t}dt
    =
    \1_{u<v}+\Ocal\biggl(
    \frac{1}{T|v-u|}+\frac{1}{Tu}+\frac{v^2u}{T^3}
    \biggr).
  \end{equation*}
\end{lem}

\begin{proof}
  The formula
  \begin{equation*}
    \frac{1}{\pi}\int_{-T}^T e^{ivt}\frac{\sin ut}{t}dt
    =
    \1_{u>v}+\Ocal\biggl(
    \frac{1}{T|v-u|}
    \biggr)
  \end{equation*}
  is classical and can for instance be found in \cite[Lemma
  2.2]{Harman*07} (and it is already used in this context by
  R.~Vaughan, see top of page~114 of \cite{Vaughan*80}). We then write
  \begin{align*}
    \frac{1}{\pi}\int_{-T}^T e^{ivt}\frac{\sin ut}{t}dt
    &=\frac{2}{\pi}\int_{0}^T (\cos(vt)-1)\frac{\sin ut}{t}dt
      +\frac{2}{\pi}\int_{0}^T \frac{\sin ut}{t}dt
    \\&=\frac{-4}{\pi}\int_{0}^T \sin^2(vt/2)\frac{\sin ut}{t}dt
      + 1 +\Ocal\biggl(\frac{1}{Tu}\biggr).
  \end{align*}
  We use $|\sin vt|\le vt$ as well as $|\sin ut|\le ut$ for the
  integral on the segment $[0,1/T]$. Since $1-1_{u>v}=1_{u<v}$, the lemma follows readily.
\end{proof}
\begin{thm}
  \label{PreciseAS}
  For any sequence $\theta$, any $T\ge1$ and any $M\ge M_0\ge 1$, we have
  \begin{multline*}
    \sum_{(n,\gP(z))=1}\mkern-12mu\theta(n)
    =\sum_{\substack{d|\gP(z)\\ d< M}}\mkern-3mu\mu(d)
    \sum_{d|n}\theta(n)
    +
    \mkern-5mu\sum_{\substack{M\le \ell\le M_0z}} \mkern-12mu\tilde{a}_\ell \theta(k\ell)
    -
    \int_{1/T}^T
    \mkern-5mu\sum_{\substack{k\ell\ge M\\ M_0\le \ell\le M_0z}}\mkern-12mu a_\ell(t)b_k(t)\theta(k\ell)
    \frac{dt}{t}
    \\
    +
    \Ocal\biggl(
    \frac{z}{T}\sum_{\substack{\ell|\gP(z)\\ M_0\le \ell<M_0 z}}
    \sum_{\substack{d|\gP(z)\\ \ell d\ge M}}\biggl|\sum_{\ell d|n}\theta(n)\biggr|
    \biggr)
  \end{multline*}
  where
  \begin{equation*}
    \tilde{a}_\ell =
    \begin{cases}
      \mu(\ell)&\text{if $\ell|\gP(z)$ and $M_0\le \ell<M_0P^{-}(\ell)$,}\\
      0&\text{otherwise,}
    \end{cases}
    ,\ 
    \mkern10mu b_k(t)=\mkern-10mu
    \sum_{\substack{p'dm=k\\ p'<z, p'\in\gP\\ d|\gP(p')}}\mu(d)
    \sin(t\log p'),
  \end{equation*}
  and
  \begin{equation*}
    a_\ell(t)=
    \begin{cases}
      \frac{4}{\pi}\mu(\ell)\sin^2\bigl(\log(P^{-}(\ell)-\tfrac12)t/2)
      &\text{if $\ell|\gP(z)$ and $M_0\le \ell<M_0P^{-}(\ell)$,}\\
      0&\text{otherwise.}
    \end{cases}
  \end{equation*}
\end{thm}
\begin{proof}
  We use Lemma~\ref{SieveBilinearSimple} with $\psi(d)=\sum_{d|n}\theta(n)$.
  We use Lemma~\ref{VFA} on $\log p'$ and $\log(P^{-}(\ell)-\frac12)$
  of Lemma~\ref{SieveBilinearSimple}
  and get
  \begin{multline*}
    \sum_{\substack{\ell|\gP(z)\\ M_0\le \ell<M_0 P^{-}(\ell)}}
    \mu(\ell)
    \sum_{\substack{p'< P^{-}(\ell)\\ p'\in\gP(z)}}
    \sum_{\substack{d|\gP(p')\\ dp'\ell\ge M}}\mu(d) 
    \psi(\ell dp')
    \\
    =
    \frac{4}{\pi}
    \int_{1/T}^T
    \sum_{\substack{\ell|\gP(z)\\ M_0\le \ell<M_0 P^{-}(\ell)}}
    \mu(\ell)\sin^2\bigl(\log(P^{-}(\ell)-\tfrac{1}{2})t/2\bigr)
    \sum_{\substack{p'<z, p'\in\gP(z)\\ d|\gP(p')\\ dp'\ell\ge M}}\mu(d)
    \sin(t\log p')
    \psi(\ell dp')
    \frac{dt}{t}
    \\
    +
    \Ocal\biggl(
    \sum_{\substack{\ell|\gP(z)\\ M_0\le \ell<M_0 P^{-}(\ell)}}
    \sum_{\substack{p'<z, p'\in\gP(z)\\ d|\gP(p')\\ dp'\ell\ge M}}|\psi(\ell dp')|
    \biggl(\frac{z}{T}+\frac{\log^3z}{T^3}\biggr)
    \biggr).
  \end{multline*}
  The theorem follows readily after the change of notation
  $dp'\rightarrow d$ in the last sum.
\end{proof}

\begin{proof}[Proof of Theorem~\ref{AS}]
  We first use Lemma~\ref{firstbase} to reduce $Z$ to  $z$ and then appeal
  to Theorem~\ref{PreciseAS}.
  \proofsubsection{First reduction}
  We simplify the error term by assuming that $\theta(n)=0$ when $n>N$
  and by first noticing that
  \begin{equation*}
    \biggl|\sum_{\ell d|n}\theta(n)\biggr|
    \le \sqrt{\frac{N}{d\ell}}\|\theta\|_2,
  \end{equation*}
  getting the error term to be $\Ocal$ of
  \begin{equation*}
    \frac{z\sqrt{N}\|\theta\|_2}{T}\sum_{\substack{\ell|\gP(z)\\ M_0\le \ell<M_0 z}}
    \sum_{\substack{d|\gP(z)\\ \ell d\le N}}\frac{1}{\sqrt{d\ell}}
    \ll
    \frac{z\sqrt{N}\|\theta\|_2}{T}\mkern-5mu
    \sum_{\substack{\ell|\gP(z)\\ M_0\le \ell<M_0
        z}}\mkern-9mu\frac{\sqrt{N/\ell}}{\sqrt{\ell}}
    \ll \frac{zN\|\theta\|_2}{T}\log(M_0z).
  \end{equation*}
  We select $T=N^3\log N$.
  \proofsubsection{Second reduction}
  On the integral over $[1/T,1]$, we apply the change of variable
  $t\mapsto 1/t$ so that we now have the sum
  \begin{equation*}
    \int_{1}^T
    \mkern-5mu\sum_{\substack{k\ell\ge M\\ M_0\le \ell\le
        M_0z}}\mkern-12mu
    a_\ell(t)b_k(t)\theta(k\ell)
    \frac{dt}{t}
    +
    \int_{1}^T
    \mkern-5mu\sum_{\substack{k\ell\ge M\\ M_0\le \ell\le
        M_0z}}\mkern-12mu
    a_\ell(1/t)b_k(1/t)\theta(k\ell)
    \frac{dt}{t}.
  \end{equation*}
  We reduce each to the interval $[0,1]$ by the change of variable
  $t\mapsto e^{y\log T}$, getting
  \begin{equation*}
     \int_0^1
    \mkern-5mu\sum_{\substack{k\ell\ge M\\ M_0\le \ell\le
        M_0z}}\mkern-12mu
    a_\ell(T^y)b_k(T^y)\theta(k\ell)
    dy
    +
    \int_{0}^1
    \mkern-5mu\sum_{\substack{k\ell\ge M\\ M_0\le \ell\le
        M_0z}}\mkern-12mu
    a_\ell(T^{-y})b_k(T^{-y})\log T\theta(k\ell)
    dy.
  \end{equation*}
  We then displace the second range of integration from 1 to 2, so
  that we may join both integrals with a proper definition of $a_\ell$
  and $b_k$. We further write
  \begin{equation*}
    \sum_{\substack{M\le \ell\le M_0z}}
    \mkern-12mu\tilde{a}_\ell \theta(k\ell)
    =
    \int_2^3 \sum_{\substack{ M\le \ell\le M_0z}}
    \mkern-12mu\tilde{a}_\ell \theta(k\ell) dy
  \end{equation*}
  so that the sum of our three bilinear contributions can be
  folded in a single integral from 0 to~3 with coefficients
  $a_\ell(t)$ and $b_k(t)$ respectively bounded in absolute value  by~$4/\pi$ and $\tau_3(k)\log T\le 4\tau_3(k)\log N$. It looks better to have an
  integration ranging from~0 to~1, something we achieve by  change
  of variable $y'=y/3$, and divide $a_\ell(t)$ by $4/\pi$, a factor we transfer in $b_k(t)$, concluding the proof.
\end{proof}
\section{Recovering the Harman Sieve}

In \cite{Harman*83,Harman*96}, G.~Harman developed his so called {\it Alternative Sieve} for primes. Here we provide another proof of a weak version of it.
\begin{thm}[Harman sieve,a  weak version]
  \label{thHS}
  Let $z\ge1$, $\mathcal{A}$ and $\mathcal{B}$ be two subsets of points on $[1,N]$ for which there exist four parameters $D\ge1$,
  $\alpha\in(0,1)$, $\lambda\in\mathbb{C}$ and $Y\ge2$ that verify the following assumptions. 
  \begin{itemize}
      \item For any sequence of complex numbers $a(d)$ with $|a(d)|\le 1$, we have
      \begin{equation*}
          \sum_{\substack{d < D, n\in\mathcal{A}\\ d|n}}
          a(d)
          =
          \lambda\sum_{\substack{d < D, n\in\mathcal{B}\\ d|n}}
          a(d)
          +\Ocal(Y).
      \end{equation*}
      \item For any two sequences of complex numbers $a(\ell)$ and $b(m)$ such that $|a(\ell)|\le 1$ and $|b(m)|\le\tau_6(m)$, we have
      \begin{equation*}
          \sum_{\substack{N^{\alpha} < \ell\le zN^{\alpha}\\ m:\ell m\in\mathcal{A}}}
          a(\ell)
          b(m)
          =
          \lambda
          \sum_{\substack{N^{\alpha} < \ell\le zN^{\alpha}\\ m:\ell m\in\mathcal{B}}}
          a(\ell)
          b(m)
          +\Ocal(Y).
      \end{equation*}
  \end{itemize}
  Then, for any sequence $(c_r)_{r\le R}$ verifying $|c_r|\le 1$ and $N^\alpha R\le D$,
  we have
  \begin{equation*}
      \sum_{r\le R}c_r
      \sum_{\substack{n\in\mathcal{A}\\ (n,\gP(z))=1\\ r|n}}1
      =
      \lambda
      \sum_{r\le R}c_r
      \sum_{\substack{n\in\mathcal{B}\\ (n,\gP(z))=1\\ r|n}}1
      +\Ocal\biggl(Y\log N+\frac{RN}{{z}}\biggr).
  \end{equation*}
\end{thm}
\noindent
We call this version ``weak" as the range in $R$ is smaller than in the full theorem of Harman, but the case $R=1$ is already of interest.
There are several versions of this theorem according to the bounds we assume on the coefficients $a$ and $b$. They are all essentially equivalent in practice. 
\begin{proof}
  Let us use Theorem~\ref{AS} with $Z=z$ and $M=M_0=N^\alpha$ twice with~$u_n$ being the characteristic function of $\mathcal{A}$ and of $\mathcal{B}$. Let us notice first that we may discard those~$r$'s which are not coprime with $\gP(z)$.
  \proofsubsection{Linear part}
  We need to bound
  \begin{align*}
     \Delta_1
     &=
      \sum_{r}c_r\sum_{\substack{d|\gP(z)\\ d<M}}\mu(d)
      \sum_{\substack{n\in\mathcal{A}\\ [d,r]|n}}1
      -
      \lambda\sum_{r}c_r\sum_{\substack{d|\gP(z)\\ d<M}}\mu(d)
      \sum_{\substack{n\in\mathcal{B}\\ [d,r]|n}}1
      \\&=
      \sum_{s<MR}c^*_s\biggl(
      \sum_{\substack{n\in\mathcal{A}\\ s|n}}1
      -
      \lambda\sum_{\substack{n\in\mathcal{B}\\ s|n}}1\biggr)
  \end{align*}
  where
  \begin{equation*}
      c_s^*=\sum_{\substack{[r,d]=s\\ d|\gP(z),d<M\\ r\le R}}c_r\mu(d).
  \end{equation*}
  Please notice that at most one decomposition $s=[r,d]$ may occur as $d|\gP(z)$ and $(r,\gP(z))=1$.
  Therefore our first hypothesis applies with $D=RM=RN^\alpha$, gives the error term~$\Ocal(Y)$, provided that $R\le D/N^\alpha$.
  \proofsubsection{Bilinear part}
  We have
  \begin{equation*}
     \Delta_2
     =
      \sum_{r}c_r
      \biggl(
      \sum_{\substack{M_0\le \ell\le M_0z\\ r|k\ell\\ k\ell\in\mathcal{A}}}
      \alpha_\ell(t)\beta_k(t)
      -
      \lambda
      \sum_{\substack{ M_0\le \ell\le M_0z\\ r|k\ell\\ k\ell\in\mathcal{B}}}
      \alpha_\ell(t)\beta_k(t)
      \biggr).
  \end{equation*}
  As $\alpha_\ell(t)$ has its support on $\ell|\gP(z)$, and $c_r$ has its support outside of this set, the condition $r|k\ell$ reduces to~$r|k$. Hence
  \begin{equation*}
     \Delta_2
     =
      \sum_{\substack{ M_0\le \ell\le M_0z\\ rm\ell\in\mathcal{A}}}
      \alpha_\ell(t)c_r\beta_{rm}(t)
      -
      \lambda
      \sum_{\substack{M_0\le \ell\le M_0z\\ rm\ell\in\mathcal{B}}}
      \alpha_\ell(t)c_r\beta_{rm}(t)
  \end{equation*}
  from which we deduce
  \begin{equation*}
     \Delta_2
     =
      \sum_{\substack{ M_0\le \ell\le M_0z\\ k\ell\in\mathcal{A}}}
      \alpha_\ell(t)\beta^*_k(t)
      -
      \lambda
      \sum_{\substack{M_0\le \ell\le M_0z\\ k\ell\in\mathcal{B}}}
      \alpha_\ell(t)\beta^*_k(t)
  \end{equation*}
  where $\beta^*_k(t)=\beta_k(t)\sum_{r|k}c_r$ obviously verifies $|\beta^*_k(t)|\ll \tau_6(k)\log N$. The theorem follows swiftly.
\end{proof}

\section{The trigonometric polynomial around $\alpha=0$}

We want to approximate $S(\beta;N)$ with enough uniformity in
$\beta$. This can be achieved by using small interval results that we
take from~\cite{Cui-Wu*14} by Z.~Cui and J.~Wu.
\begin{lem}
  \label{CuiWu}
  We have
  \begin{equation*}
    \sum_{n\le N}b(n)e(n\beta)=4C\sum_{\substack{3\le n\le N\\
        n\equiv 1[4]}}\frac{e(n\beta)}{\sqrt{\log n}}
    +\Ocal\biggl(
    \frac{N}{(\log N)^{3/2}}\bigl(1+|\beta|N^{7/10}\bigr)
    \biggr).
  \end{equation*}
\end{lem}
The initial sum carried in fact on integers that are congruent to~1 modulo~4. The constant that appears is indeed $4C$ because the sum over $n$ is
restricted to such integers.
\begin{proof}
As an application of
 Theorem~1.1 of \cite{Cui-Wu*14} with $\kappa=\alpha=1/2$, $w=0$, $\delta=1/3$ (by
Weil's bound for $\zeta(s)$), we get for $y\in[x^{7/10}/\log x,x]$ that
\begin{equation*}
  \sum_{x<n\le x+y}b(n)=\frac{Cy}{\sqrt{\log
      x}}+\Ocal\biggl(\frac{y}{(\log x)^{3/2}}\biggr).
\end{equation*}
It is not difficult to deduce from this estimate that
\begin{equation*}
  \sum_{x<n\le x+y}b(n)=4C\sum_{\substack{x<n\le x+y\\
      n\equiv1[4]}}\frac{1}{\sqrt{\log n}}
      +\Ocal\biggl(\frac{y}{(\log x)^{3/2}}\biggr).
\end{equation*}
We write, with $y=N^{7/10}/\log N$,
\begin{align*}
  \sum_{N<n\le 2N}b(n)e(n\beta)
  &=\sum_{0\le k\le N/y}\sum_{N+ky<n\le N+(k+1)y}b(n)e(n\beta)+\Ocal(y)
  \\&=\sum_{0\le k\le N/y}e(\beta(N+ky))\sum_{N+ky<n\le N+(k+1)y}b(n)
  \\&\qquad\qquad+\Ocal\biggl(\frac{|\beta|y N}{\sqrt{\log N}}+
  \frac{N}{(\log N)^{3/2}}\biggr)
\end{align*}
and after easy steps, our lemma follows.
\end{proof}

\section{Distribution in `small' arithmetic progressions via analytical means}
\label{SmallAPAM}
The distribution of Gaussian integers have been studied in several
works, starting by the main Theorem of \cite{Prachar*53} of
K.~Prachar, then Satz~1 of \cite{Rieger*65b} by G.~J.~Rieger and it is
also treated in \cite{Iwaniec*76} by H.~Iwaniec, or in Theorem~14.7 of \cite{Friedlander-Iwaniec*10} by
J.~Friedlander \& H.~Iwaniec.

We need a better error term for Theorem~\ref{Trigo}, similar to the one obtained in the paper \cite{Cui-Wu*14}. It seems this
case has not been considered, so we provide some details, though the
technique in itself is not new, see \cite{Serre*76} by J.~P.~Serre for
instance. The Selberg-Delange method as exposed
in Chapter~5, Part II of \cite{Tenenbaum*95} by G.~Tenenbaum, and more
precisely in~Theorem~3
of Section~5.3, would also lead to a solution. 
We use a path that is
intermediate and relies on \cite{Rieger*65b} by~G.~J.~Rieger.
Throughout this section, we use the following decomposition of the modulus $q\ge1$:
\begin{equation}
\label{decq}
    q=q_2q_ *,\quad \text{where $q_2$ is a power of 2 and $(q_*,2)=1$.} 
\end{equation}
Here is the main result of this section.
\begin{lem}
  \label{mainexpter}
  Let $\tau>0$ and $A\ge1$.
  Let $q=q_2q_*\le (\log N)^{\tau}$ and $a$ be prime to $q$. 
  When
  $q_*|P_{4;3}(\infty)$ and $q_2\neq4$, we have
  \begin{multline*}
    S\Bigl(\frac aq+\beta;N\Bigr)
    =\frac{4C\mu(q)}{\varphi(q)} \sum_{\substack{3\le n\le N\\
        n\equiv 1[4]}}\frac{e(n\beta)}{\sqrt{\log n}}
    \\+\Ocal\biggl(
    \frac{N2^{\omega(q)}\log q}{\varphi(q)(\log
      N)^{3/2}}\bigl(1+|\beta|N^{4/5}\bigr)\biggr)
    +\Ocal_{\tau,A}\biggl(\frac{N(1+N|\beta|)}{(\log N)^A}\biggr).
  \end{multline*}
  When
  $q_*|P_{4;3}(\infty)$ and $q_2=4$, we have
  \begin{multline*}
    S\Bigl(\frac aq+\beta;N\Bigr)
    =\frac{4Ci\chi_4(a)\mu(q_*)}{\varphi(q_*)} \sum_{\substack{3\le n\le N\\
        n\equiv 1[4]}}\frac{e(n\beta)}{\sqrt{\log n}}
    \\+\Ocal\biggl(
    \frac{N2^{\omega(q)}\log q}{\varphi(q)(\log
      N)^{3/2}}\bigl(1+|\beta|N^{4/5}\bigr)\biggr)
    +\Ocal_{\tau,A}\biggl(\frac{N(1+N|\beta|)}{(\log N)^A}\biggr).
  \end{multline*}
  When $q_*\nmid P_{4;3}(\infty)$, we have
  \begin{equation*}
    S\Bigl(\frac aq+\beta;N\Bigr)
    \ll
    \frac{N2^{\omega(q)}\log q}{\varphi(q)(\log
      N)^{3/2}}\bigl(1+|\beta|N^{4/5}\bigr)
    +\frac{N(1+N|\beta|)}{(\log N)^A}.
  \end{equation*}
    The implied constant is not effective when $\tau\ge2$ but can be
  explicitly computed otherwise.
\end{lem}
To prove this, we first prove the following lemma whose gist is that $S(\frac aq+\beta;N)$ is independent of $a$.
\begin{lem}
  \label{mainexp}
  Let $\tau>0$ and $A\ge1$.
  Let $q\le (\log N)^\tau$ and $a$ be prime to $q$. Let
  $\beta\in\mathbb{R}$. We have
  \begin{equation*}
    S\Bigl(\frac aq+\beta;N\Bigr)
    =
   \Ocal_{\tau,A}\biggl(\frac{N(1+N|\beta|)}{(\log N)^A}\biggr)
    +
    \begin{cases}
    \displaystyle
    \sum_{n\le N}\frac{c_q(n)}{\varphi(q)}b(n)e(n\beta)
    &\text{when $q_2\neq 4$},\\
    \displaystyle
    i\chi_4(a)\sum_{n\le N}\frac{c_{q_*}(n)}{\varphi(q_*)}b(n)e(n\beta)
    &\text{when $q_2=4$}. 
    \end{cases}
  \end{equation*}
  The implied constant is not effective when $\tau\ge2$ but can be
  explicitly computed otherwise.
\end{lem}
We reproduce the necessary lemmas from \cite{Rieger*65b}. In this paper,
we have $T=\exp(\sqrt{\log N})$ as written on Line~6, page~207.
There are two Dirichlet characters that behave like the principal character on our sequence $\gB$: the trivial character (and its induced versions) and the only non-trivial character modulo~$4$, say $\chi_4$, and its induced versions also denoted by~$\chi_4$, as all the elements of $\gB$ are congruent to~1 modulo~4. This character is denoted by $\chi_1$ in \cite{Rieger*65b}. Recall that, when $n$ is odd, we have $\chi_4(n)=(-1)^{\frac{n-1}{2}}$.
Here is a
combination/adaptation of Hilfssatz~11 and~12 therein.
\begin{lem}
  \label{lemRiegermult}
  There exist two positive constants $c_1$ and $c_2$ such that, for
  $N\ge 100$, $q\le \exp(\sqrt{\log N})$ and $\chi$ not induced modulo~$q$ either by the trivial character modulo~1 or by $\chi_4$, the non-principal character modulo~4, we have
  \begin{multline*}
    \sum_{n\le N}\chi(n)b(n)
    \ll
    N\exp\bigr(-c_1\sqrt{\log N}\bigr)
    \\+
    \begin{cases}
      N^{\beta(q)}(\log N)^{c_2}
      &\text{when $\chi$ is exceptional}\\
      0&\text{otherwise}.
    \end{cases}
  \end{multline*}
\end{lem}
\noindent
We can say under the hypothesis of Lemma~\ref{lemRiegermult} that for any $A>0$, we have 
\[
\sum_{n\leq N} \chi(n)b(n)
    \ll_A \frac{N}{(\log N)^A}\; .
\]
We have modified somewhat the statement of \cite{Rieger*65b}: the
modulus $k$ therein is $q$ for us. 
However we allow $q$ to be even at the price of ruling out the non-principal character modulo~4. On relying on \cite{Rieger*65b}, Hilfssatz~2 therein tells us that a polar contribution at $s=1$ occurs only in the two cases that we exclude.
We do not provide more details on  what  an
exceptional character is. We note that
$\beta(q)=1-\Ocal_\epsilon(1/q^\epsilon)$ for any positive $\epsilon$
by a famous Theorem of C.~L.~Siegel, but in an ineffective
manner. Since we shall use this theorem, some of our results are ineffective.

An important consequence emerges from this lemma: the error term is
indeed very small when we consider the non-principal character. We now
have to consider the contribution of the principal character, but we
see that its contribution is going to be independent of~$a$,
for~\eqref{baseRiegerstar}. As the
principal character is indeed real valued, we can now rely on
\cite{Serre*76} or on the Selberg-Delange Method.

\begin{lem}
  \label{GS}
  We have $\tau_q(\chi_0)=\mu(q)$, and, when $4|q$, we have $\tau_q(\chi_4)=-2i\mu(q/2)$, where
  $\tau_q(\chi)=\sum_{b\mode q}\chi(b)e(b/q)$.
\end{lem}

\begin{proof}
  The value of the Gauss sum on the trivial character is classical. Concerning $\chi_4$, we split the modulus $q$ in $q=q_2q'$ where $q'$ is odd and $q_2$ is a power of~2.
  We have $\tau_q(\chi_4)=\tau_{q_2}(\chi_4)\tau_{q'}(\chi_0)$. We readily check that $\tau_{q_2}(\chi_4)=0$ when $8|q_2$. Finally, we see that $\tau_{4}(\chi_4)=e(1/4)-e(3/4)=2i$. We may fold these results in the formula $\tau_q(\chi_4)=-2i\mu(q/2)$, concluding the proof of this lemma.
\end{proof}

\begin{proof}[Proof of Lemma~\ref{mainexp}]
We set $\alpha=\frac aq+\beta$. 
Let $q \ge 3$ be a modulus and $a$ be prime to~$q$. Let us define
  \begin{equation}
    \label{eq:4}
    S_q(\alpha; N)=\sum_{\substack{n\le N\\ (n,q)=1}}b(n)e(n\alpha).
  \end{equation}
  The reader should notice that $S_q$ depends only on the prime factors of~$q$ and not on the actual value of this parameter.
  With the aim of evaluating this function, we replace the
  additive character with a linear combination of multiplicative
  ones. To do this, we first notice that, for $n$ prime to~$q$, we have:
  \begin{align*}
    e(na/q)
    &=\sum_{c\mod^*q}\frac{1}{\varphi(q)}
    \sum_{\chi\mod q}\overline{\chi}(c)\chi(n) e(ca/q)
    \\&=
    \frac{1}{\varphi(q)}\sum_{\chi\mod q}
    \sum_{c\mod^*q}\overline{\chi}(c)e(ca/q)\chi(n)
    \\&=\frac{1}{\varphi(q)}\sum_{\chi\mod q}\tau_q(\overline{\chi})\chi(a)\chi(n).
  \end{align*}
  Therefore
  \begin{equation}
    \label{baseRiegerstar}
    S_q(a/q; N)=\frac{1}{\varphi(q)}\sum_{\chi\mod
      q}\tau_q(\overline{\chi})\chi(a)
    \sum_{n\le N}b(n)\chi(n).
  \end{equation}
  Please note that each Gauss sum $\tau_q(\overline{\chi})$ is of
  modulus at most $\sqrt{q}$.

\vspace{2mm}
\noindent
By Lemma~\ref{lemRiegermult} with $A+\tau/2$ rather than $A$, when $4\nmid q$, we have
\begin{equation}
  \label{eq:5}
  S_q(a/q; N)-\frac{\mu(q)}{\varphi(q)}S_q(0;N)
  \ll_{\tau,A} \frac{N}{(\log N)^{A}}
\end{equation}
valid for any $A\ge0$, any $q\le (\log N)^\tau$ and any $a$ prime to $q$.

When $4|q$, then $\mu(q)=0$ but the character $\chi_4$ has some effect and recalling Lemma~\ref{GS}, we get
\begin{equation}
  \label{eq:5+}
  S_q(a/q; N)+\frac{2i\mu(q/2)\chi_4(a)}{\varphi(q)}S_q(0;N)
  \ll_{\tau,A} \frac{N}{(\log N)^{A}}.
\end{equation}
Let us set
\begin{equation}
    s(q, a)=\begin{cases}
      \mu(q)&\text{when $q_2\neq 4$}\\
      -2i\chi_4(a)\mu(q/2)&\text{when $q_2=4$}.
    \end{cases}
\end{equation}

\proofsubsection{First extension}
We need to extend this inequality to $S_{qq'}(a/q; N)$ for any integer
$q'$ prime to~$q$.
It is better for later usage to change the notation.
For integers $u$ and $v$ with $(u,v)=1$, we have 
\begin{align}
  S_{uv}(f/u; N)
  &=\notag
  \sum_{w|v}\mu(w)b(w)S_{u}(fw/u; N/w)
  \\&=\notag
  \frac{1}{\varphi(u)}\sum_{w|v}\mu(w)b(w)
    s(u,fw)
  S_{u}(0;
  N/w)
  +\Ocal_{A,\tau}\biggl((\log v)\frac{N}{(\log N)^A}\biggr)
  \\&=\label{FE}
  M_1(N,u, uv,f\mod4)
  +
  \Ocal_{A,\tau}\biggl((\log v)\frac{N}{(\log N)^A}\biggr)
\end{align}
for some function $M_1$.
When $u_2\neq 4$, we readily find that $M_1(N,u, uv, f\mod 4)=\mu(u)S_{uv}(0;N)/\varphi(u)$, while, when $u_2=4$, the expression is a bit more intricate, but still only depends on $f\mod 4$.
\proofsubsection{Second extension}
  We now have to dispense with the coprimality condition $(n,q)=1$. To
  do so, we write $\delta|q^\infty$ to say that the prime factors of
  $\delta$ divide~$q$. It is possible to write every integer $n$ in
  the form $n=\delta m$ with $\delta|q^\infty$ and $(m,q)=1$, and this
  can be done in a unique manner. Since the function $d\mapsto b(d)$ is completely multiplicative, we find that
  \begin{equation*}
      S(a/q; N)
    =\sum_{\substack{\delta|q^\infty\\ \delta\le N/3}}
    b(\delta) S_q(a\delta /q; N/\delta).
  \end{equation*}
  Notice that we may assume $\delta$ to be odd. We have
  \begin{equation*}
      S_q(a\delta /q; N/\delta)
      =S_q\biggl(\frac{a\delta/(q,\delta) }{q/(q,\delta)};N/\delta\biggr)
      =S_{\frac{q}{(q,\delta)}r_q(\delta)}\biggl(\frac{a\delta/(q,\delta) }{q/(q,\delta)};N/\delta\biggr)
  \end{equation*}
  where $r_q(\delta)=\prod_{p|q,p\nmid q/(q,\delta)}p$. Notice that $r_q(\delta)$ may be smaller than $(q,\delta)$, which is why we need this extra notation. Let us resume the main evaluation:
  \begin{align}
    S(a/q; N)
    &=\sum_{\substack{\delta|q^\infty\\ \delta\le N/3}}
    b(\delta) S_{\frac{q}{(q,\delta)}r_q(\delta)}\biggl(\frac{a\delta/(q,\delta) }{q/(q,\delta)};N/\delta\biggr)\notag
    \\&\label{firstone}
    =M_2(q,N)
    +\Ocal\biggl(N
    \sum_{\substack{\delta|q^\infty\\ \delta\le N/3}}
    \frac{1}{\delta} \frac{1}{(\log(N/\delta))^A}\biggr)
  \end{align}
  where $M_2(q,N)$ is the 'main term' and is defined as
  \begin{equation}
    \label{defM}
    M_2(q,N)
    =
    \sum_{\substack{\delta|q^\infty\\ \delta\le N/3\\ (\delta,2)=1}}
    b(\delta) M_1\biggl(N, \frac{q}{(q,\delta)}, 
    \frac{q}{(q,\delta)}r_q(\delta), \frac{a\delta}{(q,\delta)}\mod 4\biggr)
  \end{equation}
  by appealing to~\eqref{FE}. We find again that, when $q_2\neq4$, this main term is independent of~$a$, while, when $q_2=4$, it depends only on the class of $a$ modulo~4. 
  Concerning the error term in~\eqref{firstone},
  when $\delta\ge \sqrt{N}$, we use
  \begin{align*}
    \sum_{\substack{\delta|q^\infty\\ \sqrt{N}\le \delta\le N/3}}
    \frac{1}{\delta}
    \frac{1}{(\log(N/\delta))^A}
    &\le \sum_{\substack{\delta|q^\infty}}
    \frac{1}{\delta}\biggl(\frac{\delta}{\sqrt{N}}\biggr)^{1/2}
    \\&\le \frac{1}{N^{1/4}}(1-1/\sqrt{2})^{-\omega(q)}
    \ll
    \frac{1}{N^{1/4}} 4^{\frac{\log q}{\log 2}}
    \\&\ll \frac{1}{N^{1/4}} 8^{\tau \log\log N}
    \ll_{\tau,A} \frac{1}{(\log N)^A}
  \end{align*}
  since $q\le (\log N)^\tau$. For the smaller $\delta$, we use
  \begin{equation*}
    \sum_{\substack{\delta|q^\infty\\  \delta\le \sqrt{N}}}
    \frac{1}{\delta}
    \frac{1}{\log(N/\delta)^A}
    \ll
      \frac{1}{(\log N)^A}\prod_{p|q}(1+1/p)
    \ll \frac{\log q}{(\log N)^A}.
  \end{equation*}
  We conclude that
  \begin{equation*}
    S(a/q;N)-M_2(q,N)
    \ll_{\tau,A}
    \frac{N\log q}{(\log N)^A}.
  \end{equation*}
  The expression for $M_2(q,N)$ is rather unwieldy but discloses that this quantity does not depend on~$a$. We may therefore sum over $a'$ congruent to $a$ modulo $q/q_*$, and
  prime to $q_*$. On changing $A$ to
  $A+\tau+1$, we find that
  \begin{equation}
    \label{finalasurq}
    S(a/q;N)=
    \Ocal_{\tau,A}\biggl(\frac{N}{(\log N)^A}\biggr)
    +
    \begin{cases}
    \sum_{n\le N}\frac{c_q(n)}{\varphi(q)}b(n)
    &\text{when $q_2\neq4$},\\
    i\chi_4(a)\sum_{n\le N}\frac{c_{q_*}(n)}{\varphi(q_*)}b(n)
    &\text{when $q_2=4$}.
    \end{cases}
  \end{equation}
\proofsubsection{Analytic variations}
  
  We readily get: 
  \begin{align*}
    S(\alpha;N)
    &=
      \sum_{n\le N}b(n)e(na/q)e(n\beta)
    \\&=
    \sum_{n\le N}b(n)e(na/q)\biggl(2i\beta\pi\int_0^{n}e(\beta t)dt+1\biggr)
    \\&=
    2i\pi\beta\int_0^{N}e(\beta
    t) (S(a/q;N)-S(a/q;t))dt+S(a/q;N).
  \end{align*}
  On using~\eqref{finalasurq}, $c_q(n)=\sum_{d|n,d|q}d\mu(q/d)$ and the fact that the $d\mapsto b(d)$ function is completely multiplicative, we find that
  \begin{multline*}
    S(\alpha;N)
    =
      \sum_{\substack{d|q}}\frac{\mu(q/d)db(d)}{\varphi(q)}\biggl(
      2i\pi\beta\int_0^{N}e(\beta
      t) \Bigl(B\Bigl(\frac{N}{d}\Bigr)-B\Bigl(\frac{t}{d}\Bigr)dt+B\Bigl(\frac{N}{d}\Bigr)
      \biggr)
     \\  +\Ocal_{\tau,A}\biggl(\frac{N(1+N|\beta|)}{(\log N)^A}\biggr)
  \end{multline*}
  where $B(t)=\sum_{n\leq t}b(n)$ is defined in~\eqref{refC}.
  The integral may also be written in the form
  \begin{equation*}
    \int_0^{N/d}e(d\beta
      u) (B(N/d)-B(u))du+B(N/d).
  \end{equation*}
  Lemma~\ref{mainexp} follows swiftly.
\end{proof}
We now first restate Lemma~\ref{mainexp} in the format used in Lemma~\ref{CuiWu}.
\begin{lem}
  \label{mainexpbis}
  Let $\tau>0$ and $A\ge1$.
  Let $q\le (\log N)^{\tau}$, $a$ prime to $q$ and let $q^*$ be the largest
  divisor of~$q$ such that $b(q^*)=1$. When $q_2\neq 4$, we have
  \begin{multline*}
    S\Bigl(\frac aq+\beta;N\Bigr)
    =4C  \sum_{\substack{3\le n\le N\\
        n\equiv 1[4]}}\frac{e(n\beta)}{\sqrt{\log n}}
    \frac{\mu(q/q^*)c_{q^*}(n)}{\varphi(q)}
    \\+\Ocal\biggl(
    \frac{N2^{\omega(q)}\log q}{\varphi(q)(\log
      N)^{3/2}}\bigl(1+|\beta|N^{4/5}\bigr)\biggr)
    +\Ocal_{\tau,A}\biggl(\frac{N(1+N|\beta|)}{(\log N)^A}\biggr).
  \end{multline*}
  When $q_4= 4$, we have
  \begin{multline*}
    S\Bigl(\frac aq+\beta;N\Bigr)
    =4C  i\chi_4(a)\sum_{\substack{3\le n\le N\\
        n\equiv 1[4]}}\frac{e(n\beta)}{\sqrt{\log n}}
    \frac{\mu(q_*/q^*)c_{q^*}(n)}{\varphi(q_*)}
    \\+\Ocal\biggl(
    \frac{N2^{\omega(q)}\log q}{\varphi(q)(\log
      N)^{3/2}}\bigl(1+|\beta|N^{4/5}\bigr)\biggr)
    +\Ocal_{\tau,A}\biggl(\frac{N(1+N|\beta|)}{(\log N)^A}\biggr).
  \end{multline*}
\end{lem}

\begin{proof}
  Let us set, for this proof,
  \begin{equation*}
     M(q,N,\beta)=
    \sum_{\substack{d|q}}\frac{\mu(q/d)db(d)}{\varphi(q)}S(d\beta;
    N/d).
  \end{equation*}
  This is the main term obtained in Lemma~\ref{mainexp} when $q_2\neq4$. Otherwise the main term is $\chi_4(a)M(q_*,N,\beta)$.
  By Lemma~\ref{CuiWu}, we find that, when $q\le N^{1/10}$, we have
  \begin{equation*}
    M(q,N,\beta)=
      4C \sum_{\substack{d|q}}\frac{\mu(q/d)db(d)}{\varphi(q)}
      \sum_{\substack{3\le n\le N/d\\
    n\equiv 1[4]}}\frac{e(nd\beta)}{\sqrt{\log n}}
    +\Ocal\biggl(
    \frac{N2^{\omega(q)}(1+|\beta|N^{4/5})}{\varphi(q)(\log N)^{3/2}}
    \biggr).
  \end{equation*}
  To replace $\log n$ by $\log (nd)$, we note that
  \begin{equation*}
    \frac{1}{\sqrt{\log n}}-\frac{1}{\sqrt{\log (nd)}}
    =\frac{\log d}{\sqrt{(\log n)(\log(nd))}(\sqrt{\log n}+\sqrt{\log
        (nd)})}
    \le \frac{\log d}{(\log n)^{3/2}},
  \end{equation*}
  from which we deduce that
  \begin{align*}
    M(q,N,\beta)
    &=
    4C  \sum_{\substack{3\le n\le N\\
        n\equiv 1[4]}}\frac{e(n\beta)}{\sqrt{\log n}}
    \sum_{\substack{d|q\\ d|n}}\frac{\mu(q/d)db(d)}{\varphi(q)}
    +\Ocal\biggl(
    \frac{N2^{\omega(q)}\log q}{\varphi(q)(\log
    N)^{3/2}}
    \bigl(1+|\beta|N^{4/5}
    \bigr)\biggr)
    \\&=
    4C  \sum_{\substack{3\le n\le N\\
        n\equiv 1[4]}}\frac{e(n\beta)}{\sqrt{\log n}}
    \frac{\mu(q/q^*)c_{q^*}(n)}{\varphi(q)}
    +\Ocal\biggl(
    \frac{N2^{\omega(q)}\log q}{\varphi(q)(\log N)^{3/2}}\bigl(1+|\beta|N^{4/5}\bigr)\biggr).
  \end{align*}
  as announced.
\end{proof}


\begin{proof}[Proof of Lemma~\ref{mainexpter}]
  Our main remark is the following estimate, valid provided $d\nmid2$:
  \begin{equation*}
    \sum_{\substack{n\le M\\ n\equiv 1[4]}}c_d(n)
        =
        \sum_{\substack{n\le M\\ n\equiv 1[4]}}\sum_{\delta|n,\delta|d}\mu(d/\delta)\delta
         =
         \sum_{\delta|d}\delta\mu(d/\delta)
        \sum_{\substack{n\le M\\ n\equiv 1[4]\\ \delta|n}}1
   \ll \phi_+(d)
  \end{equation*}
  where
  \begin{equation}
    \label{defphiplus}
    \phi_+(d)=d\prod_{p|d}(1+1/p).
  \end{equation}
  We deduce for this bound that $\sum_{\substack{n\le M\\ n\equiv
      1[2]}}e(n\beta)c_d(n)\ll (1+|\beta| N)\phi_+(d)$. Hence, the
  main term given in Lemma~\ref{mainexpbis} is
  small when $q^*\neq1$.
\end{proof}

\section{A local model}
\label{Local}
We consider some (large) positive real number $K$ and put
\begin{equation}
  \label{defSflatK}
  S^\flat_K(\alpha;N)=
  4C\prod_{\substack{p\le K\\ p\equiv 3[4]}}\biggl(1-\frac1p\biggr)^{-1}
  \sum_{\substack{3\le n\le N\\ n\equiv1[4]\\ (n,P_{4;3}(K))=1}}
  \frac{e(n\alpha)}{\sqrt{\log n}}.
\end{equation}
\begin{lem}
  \label{easycase}
  When $4\alpha\notin\mathbb{Z}$, we have
  $\displaystyle
  \sum_{\substack{3\le n\le N\\ n\equiv1[4]}}
  \frac{e(\alpha n)}{\sqrt{\log n}}\ll 1/\|4\alpha\|$.
\end{lem}

\begin{proof}
  Indeed, by summation by parts, we find that the sum, say
  $\Sigma(\alpha)$, may be rewritten as
  \begin{equation*}
    \Sigma(\alpha)
    =
    \int_3^N\sum_{\substack{3\le n\le t\\ n\equiv1[4]}}
    e(\alpha n)\frac{dt}{t(\log t)^{3/2}}
    +
    \frac{1}{\sqrt{\log N}}\sum_{\substack{3\le n\le N\\ n\equiv1[4]}}
      e(\alpha n)
    \ll1/\|4\alpha\|.
  \end{equation*}
\end{proof}

\begin{lem}
  \label{helper}
  When $\alpha=\frac{a}{q}+\beta$, $q|4d$ and $d|P_{4;3}(\infty)$, let us define
  \begin{equation*}
      \Sigma_d(N;\alpha)
      =\sum_{\substack{3\le n\le N\\ n\equiv1[4]\\ d|n}}\frac{e(n\alpha)}{\sqrt{\log
    n}}.
  \end{equation*}
  We have $\Sigma_d(N;\alpha)=e(aq_*/q_2)\Sigma_q(N;\beta)$
  which is therefore independent on $a$ prime to $q$ when $q_2\neq 4$ and depends only on $a\mod 4$ when $q_2=4$.
\end{lem}

\begin{proof}
  Let $n_0\ge 3$ be the smallest integer solution to $n\equiv1[4]$ and $d|n$.  
  Let us write
  \begin{equation*}
      \frac{a}{q}\equiv\frac{a_2}{q_2}+\frac{a_*}{q_*}\ \mod 1\; .
  \end{equation*}
  We have
  \begin{align*}
      \frac{na}{q}
      &\equiv\frac{n_0a}{q}+\frac{(n-n_0)a}{q}
      \equiv \frac{n_0a}{q}\ \mod 1
      \\& \equiv\frac{n_0a_2}{q_2}+\frac{n_0 a_*}{q_*} \equiv\frac{n_0a_2}{q_2}\ 
      \mod 1
  \end{align*}
  because $q_*|d|n_0$. We also know that $q_2|4$ (as $q|4d$), that $n_0\equiv 1[4]$ and that $a_2\equiv a q_*^{-1}[q_2]\equiv a q_*[q_2]$ because, modulo~2 or~4, a class equals its inverse class.
\end{proof}

\begin{lem}
  \label{model}
  Let $Q=P_{4;3}(K)$. Let $\alpha\in\mathbb{R}/\mathbb{Z}$ be given with a rational approximation $\alpha-a/q=\beta$, where $q$ is a positive integer, $a$ is coprime with $q$ and $qQ|\beta|\le 1$. Write $q=q_2q_ *$ with $q_2$ a power of $2$ and $q_*$ odd. We have, when
  $q_2\nmid4$ or $q_ *\nmid P_{4;3}(K)$:
  \begin{equation*}
    S^\flat_K(\alpha;N)\ll qe^{K}.
  \end{equation*}
  When $q_2\neq4$ and $q_*\mid P_{4;3}(K)$, we find that
  \begin{equation*}
    S^\flat_K(\alpha;N)=\frac{4\mu(q)C}{\varphi(q)}
    \sum_{\substack{3\le n\le N\\
        n\equiv1[4]}}\frac{e(n\beta)}{\sqrt{\log n}}
    +\Ocal\biggl(
    \frac{e^{-K/4}N }{q\sqrt{\log N}}+qe^{K/4}\biggr).
  \end{equation*}
  When $q_2=4$ and $q_*\mid P_{4;3}(K)$, we find that
  \begin{equation*}
    S^\flat_K(\alpha;N)=\frac{4\mu(q_*)Ci\chi_4(a)}{\varphi(q_*)}
    \sum_{\substack{3\le n\le N\\
        n\equiv1[4]}}\frac{e(n\beta)}{\sqrt{\log n}}
    +\Ocal\biggl(
    \frac{e^{-K/4}N }{q\sqrt{\log N}}+qe^{K/4}\biggr).
  \end{equation*}
  
\end{lem}

\begin{proof}
  We set $\delta(K)=\varphi(P_{4;3}(K))/P_{4;3}(K)$.
  We readily find that
  \begin{align*}
    \delta(K)\frac{S^\flat_K(\alpha;N)}{4C}
    &=
      \sum_{d|P_{4;3}(K)}\mu(d)
      \sum_{\substack{3\le n\le N\\ n\equiv1[4]\\ d|n}}\frac{e(n\alpha)}{\sqrt{\log n}}
    \\&=
      \sum_{\substack{d|P_{4;3}(K)\\ q|4d}}\mu(d)
      \sum_{\substack{3\le n\le N\\ n\equiv1[4]\\ d|n}}\frac{e(n\alpha)}{\sqrt{\log
    n}}
    +\Ocal\Bigl(\sum_{d|P_{4;3}(K)}q\Bigr)
    \end{align*}
  by Lemma~\ref{easycase}.
  Notice that, by the prime number theorem modulo~4, we have
  $\sum_{d|P_{4;3}(K)}1\ll e^{(1+o(1))(\log 2)\frac{K}{2\log K}}\ll_\ve e^{\ve K}$ for any $\ve >0$. Let us write $q=q_2q_*$ where $q_2$ is a power of 2 and $q_*$ is odd.
  The first sum of course is empty when $q_*\nmid
  P_{4;3}(K)$ or when $q_2\nmid 4$. 
  Lemma~\ref{helper} 
  implies that the main term above does not depend on~$a$ when $q_2\neq4$, and only on $a\mod 4$ when $q_2=4$. In the first case, we may thus average over all~$a$ prime to~$q$, and over all~$a'$ in $(\mathbb{Z}/4q_*\mathbb{Z})^\times$ congruent to~$a$ modulo~4.
  In the first case, this gives us, when $q_*|P_{4;3}(K)$:
  \begin{equation*}
    \delta(K)\frac{S^\flat_K(\alpha;N)}{4C}
    =
    \sum_{\substack{3\le n\le N\\ n\equiv1[4]\\ (n,P_{4;3}(K))=1}}
    \frac{e(n\beta)c_q(n)}{\varphi(q)\sqrt{\log n}}
    +\Ocal_\ve(qe^{\ve K})
    .
  \end{equation*}
  But as $q$ and $n$ are coprime, the dependence in $q$ can simply be
  factored out, getting
  \begin{equation*}
    \delta(K)\frac{S^\flat_K(\alpha;N)}{4C}
    =
    \frac{\mu(q)}{\varphi(q)}
    \sum_{\substack{3\le n\le N\\ n\equiv1[4]\\ (n,P_{4;3}(K))=1}}
    \frac{e(n\beta)}{\sqrt{\log n}}
    +\Ocal(qe^{\ve K}).
  \end{equation*}
  To further study the main term, we assume that $\mu(q)\neq0$.
  We use summation by parts to
  infer that
  \begin{align*}
    \sum_{\substack{3\le n\le N\\ n\equiv1[4]\\ (n,P_{4;3}(K))=1}}
    \frac{e(n\beta)}{\sqrt{\log n}}
    &=
    \int_{3}^N \sum_{\substack{3\le n\le t\\ n\equiv1[4]\\ (n,P_{4;3}(K))=1}}1
    \biggl(2i\pi\beta \frac{e(t\beta)}{\sqrt{\log t}}+ \frac{e(\beta t)}{2t(\log t)^{3/2}}\biggr) dt
    \\&\quad +
    \sum_{\substack{3\le n\le N\\ n\equiv1[4]\\ (n,P_{4;3}(K))=1}}1\frac{e(\beta N)}{\sqrt{\log N}}.
  \end{align*}
  We employ the approximation
  $\sum_{\substack{3\le n\le t, n\equiv1[4]\\ (n,P_{4;3}(K))=1}}1=\delta(K)\frac{t}{4}+\Ocal_\ve(e^{\ve K})$ and reach
  \begin{align*}
    \sum_{\substack{3\le n\le N\\ n\equiv1[4]\\ (n,P_{4;3}(K))=1}}
    \frac{e(n\beta)}{\sqrt{\log n}}
    &=
    \delta(K)\int_{3}^N \sum_{\substack{3\le n\le t\\ n\equiv1[4]}}1
    \biggl(2i\pi\beta \frac{e(t\beta)}{\sqrt{\log t}}+ \frac{e(\beta t)}{2t(\log t)^{3/2}}\biggr) dt
    \\&\quad +
    \delta(K)\sum_{\substack{3\le n\le N\\ n\equiv1[4]}}1
    \frac{e(\beta N)}{\sqrt{\log N}}
    +\Ocal_\ve
    \biggl(
    \frac{e^{\ve K}(1+N|\beta|)}{\sqrt{\log N}}\biggr)
    .
  \end{align*}
  The approximation is the same quantity with $K=1$, therefore 
  \begin{equation*}
    {\mu(q)}{\varphi(q)}\delta(K)\frac{S^\flat_K(\alpha;N)}{4C}
    =\delta(K)
    \sum_{\substack{3\le n\le N\\ n\equiv1[4]}}
    \frac{e(n\beta)}{\sqrt{\log n}}
    +\Ocal_\ve \biggl(
    \frac{e^{\ve K}|\beta|N}{\sqrt{\log N}}+qe^{\ve K}
    \biggr).
  \end{equation*}
  Finally, we notice that $\delta(K)^{-1}\ll \log K$ and that $|\beta|\ll e^{-K/3}/q$. The lemma is proved in the case $q_2\neq4$.
  
  When $q_2=4$, 
  we skip some steps and reach
  \begin{equation*}
      \delta(K)\frac{S^\flat_K(\alpha;N)}{4C}
    =
    i\chi_4(a)\frac{\mu(q_*)}{\varphi(q_*)}
    \sum_{\substack{3\le n\le N\\ n\equiv1[4]\\ (n,P_{4;3}(K))=1}}
    \frac{e(n\beta)}{\sqrt{\log n}}
    +\Ocal(qe^{\ve K}).
  \end{equation*}
  The lemma follows readily. 
\end{proof}

\begin{lem}\label{L1b}
  For any $\ell>2$ and $K\ge2$, we have
  \begin{equation*}
    \int_0^1\biggl|S^\flat_K(\alpha;N)  \biggr|^\ell d\alpha\ll
    \sqrt{\log K}\frac{N^{\ell-1}}{(\log N)^{\ell/2}}.
  \end{equation*}
\end{lem}
Anticipating on Lemma~\ref{L1}, the above looses a factor $\sqrt{\log K}$.

\begin{proof}
  We simply use Parseval:
  \begin{align*}
    \int_0^1\biggl|
    S^\flat_K(\alpha;N)|^\ell d\alpha
    &\le \bigl(S^\flat_K(0;N)\bigr)^{\ell-2}
    C^2\prod_{\substack{p\le K\\ p\equiv 3[4]}}\biggl(1-\frac1p\biggr)^{-2}
  \sum_{\substack{3\le n\le N\\ n\equiv1[2]\\ (n,P_{4;3}(K))=1}}
    \frac{1}{\log n},
    \\&\ll
    \delta(K)^{-\ell}\frac{(\delta(K)N)^{\ell-2}}{(\log
    N)^{\frac{\ell-2}{2}}}\delta(K)\frac{N}{\log N}
  \end{align*}
  where $\delta(K)=\prod_{\substack{p\le K\\ p\equiv
      3[4]}}(1-1/p)$. The lemma follows readily.
\end{proof}

\section{The trigonometric polynomial in the large}
Here are Lemma 3 and a modification of Lemma 2 of \cite{Vaughan*77c} by R.~Vaughan.
\begin{lem}
  \label{VaughanLem3}
  Assume that $|q\alpha-a|\le 1/q$ for some $a$ prime to~$q$. Then we have
  \begin{equation*}
      \sum_{u\le U}\max_{Z\le UV/u}
      \biggl|\sum_{v\le Z}a_ue(uv\alpha)\biggr|
      \le \|a\|_\infty \biggl(\frac{UV}{q}+U+q\biggr)
      \log(UVq).
  \end{equation*}
\end{lem}
\begin{lem}
  \label{VaughanLem2}
  Assume that $|q\alpha-a|\le 1/q$ for some $a$ prime to~$q$. Then we have
  \begin{equation*}
      \sum_{u\le U}\max_{Z\le Y}
      \biggl|\sum_{v\le Z}a_u b_v e(uv\alpha)\biggr|
      \le \|a\|_2\|b\|_ 2\log(UYq)^{5/2}
      \sqrt{UYq^{-1}+U+Y+q}
     .
  \end{equation*}
\end{lem}

It is convenient to extend these two lemmas in the form given in Lemma~4 in
\cite{Harman*83} by G.~Harman. The proof of which lies in fact between
Eq.~(15) and~(16) of \cite{Vaughan*77c}. We are slightly more precise
than G.~Harman in the dependence in the coefficients.

\begin{lem}
  \label{HarmanLem4}
  Assume that $|q\alpha-a|\le 1/q$ for some $a$ prime to~$q$. Then for $R\ge 1/2$, $U^\prime\leq 2U, V^{\prime}\leq 2V$ and $X^{\prime}\leq 2X$ we have
  \begin{multline*}
      \sum_{R<r\le 2R}
      \biggr|\sum_{\substack{U<u\le U'\\ V<v\le V'\\X<uv\le X'}}
      a_u b_v e(ruv\alpha)\biggr|
      \\
      \le \|a\|_{2,2}\|b\|_ 2R\sqrt{X}
      \biggl(
      \frac{1}{\sqrt{q}}+\frac{\sqrt{q}}{\sqrt{RX}}
      +
      \frac{1}{\sqrt{V}}+
      \frac{1}{\sqrt{RU}}
      \biggr)
      \log(RXq)^{3}
    \end{multline*}
    where the norm $\|a\|_{2,2}$ is defined by $\|a\|_{2,2}^2=\sum_{u}d(u)|a_u|^2$,
    the function $d(u)$ being the number of divisors of~$u$.
\end{lem}
\noindent
We state this inequality with a symmetric left-hand side:
$U$ and $V$ may thus be interchanged to decrease the
right-hand side.
The size conditions $u\in(U,U']$ and $v\in(V,V']$ can
easily be carried by altering the sequences $(a_u)$~and $(b_v)$. It is better to use the above when $U\le
X^{1/2}$ and to reverse the roles of $U$ and $V$ otherwise.
\begin{proof}
  Let us denote by $S$ the sum to be studied. We write
  \begin{align*}
      S 
      &= \sum_{R<r\le 2R}\sum_{U<u\le 2U} a_u 
      \sum_{rX<ruv\le rX'}
     c_r b_v \1_{V<v\le V'} e(ruv\alpha)
      \\&\le 2\sum_{UR< w \le 4UR} A_w \max_{Z\le 2RX}
      \biggl|\sum_{v\le Z/w} b^*_v  e(wv\alpha)\biggr|
  \end{align*}
  with an obvious notation.
  
 Now we can apply  Lemma~\ref{VaughanLem2}. We first notice that
  \begin{align*}
      \|A\|_2^2
      &\le \sum_{w}\Bigl(\sum_{ur=w}|a_u|\Bigr)^2
      \le \sum_{u}|a_u|^2\sum_{u|w}d(w)
      \\&\ll RU\log(RU)\sum_{u}\frac{d(u)|a_u|^2}{u}
      \ll R\sum_{u}d(u)|a_u|^2\log(RU).
  \end{align*}
  The lemma follows readily.
\end{proof}

The aim of this section is the following crucial lemma.
\begin{lem}
  \label{TrigoLargeII}
  Let $|q\alpha-a|\le 1/q$ with $q\le N$ and $(a,q)=1$. We have
  \begin{equation*}
    \frac{S(\alpha;N)}{N/\sqrt{\log N}}\ll
    \biggl(
    \frac{1}{\sqrt{q}}
    +
    \frac{1}{N^{1/6}}
      +
    \sqrt{\frac{q}{N}}
    \biggr)(\log N)^{7}
    .
  \end{equation*}
\end{lem}
As the same method would work for the primes with the same result, it
is interesting to compare with \cite[Theorem 1]{Vaughan*77c} of
R.~Vaughan. If we forget the $\log$-powers, our result has the
worse $N^{-1/6}$ compared with $N^{-1/5}$ there.

\begin{proof}
We start by localizing the variable $n$ and set
$S^*(\alpha,N)=\sum_{N/2<n\le N}b(n)e(n\alpha)$.
Theorem~\ref{AS} with $Z=\sqrt{N}$ proposes a decomposition of $S^*$ into
two types of sum
that we study one after the other. We select $M=M_0$ and set
$\alpha=\frac{a}{q}+\frac{\beta}{q^2}$.
\proofsubsection{Study of the linear part}
Lemma~\ref{VaughanLem3} gives us directly
\begin{equation}
  \label{est1}
  \sum_{\substack{d|P_{4;3}(z)\\
      d\le M}}\mu(d)\sum_{\substack{n\equiv0[d]\\ n\equiv 1[4]\\ N/2<n\le N}}e(n\alpha)
  \ll \biggl(\frac{N}{q}+M+q\biggr)\log(Nq).
\end{equation}

\proofsubsection{Study of the bilinear part}

We localize the variable $\ell$ in a diadic interval
$[L,L']\subset [M,Mz]$ with $L'\le 2L$. Let $a^*_\ell$ denotes either
$\tilde{a}_\ell$ or $a_\ell(t)$ for some $t$, and $b^*_k=b_k(t)$
for some $t$. 
Lemma~\ref{HarmanLem4} applies with $R=1/2$ and yields the upper bound
\begin{equation}
  \label{est3}
  N(\log qN)^{9/2}\biggl(\frac{1}{\sqrt{q}}
  +\frac{\sqrt{q}}{\sqrt{N}}
  +\frac{1}{\sqrt{\min(L, N/L)}}\biggr)\; .
\end{equation}

\proofsubsection{Gathering the estimates}
On gathering~\eqref{est1} and \eqref{est3} together with
the error term given by Theorem~\ref{AS}, we reach
\begin{multline*}
  S^*(\alpha,N)/N\ll
  \biggl(\frac{1}{\sqrt{q}}
  +\frac{\sqrt{q}}{\sqrt{N}}
  +\frac{1}{\sqrt{z}}+\frac{1}{\sqrt{M}}+\frac{\sqrt{Mz}}{\sqrt{N}}\biggr)
  (\log N)^{11/2}\log T
  \\+
  \frac{M\log N}{N}
  +\frac{\log N}{q}
  +\frac{z(\log N)^3}{T}.
\end{multline*}
We select $M=z=N^{1/3}$ and $T=N$, getting
\begin{equation*}
  S^*(\alpha,N)\ll
  \biggl(\frac{N}{\sqrt{q}}
  +N^{5/6}
  + \sqrt{Nq}\biggr)(\log N)^{13/2}.
\end{equation*}
We sum diadically over $N$ to
get the claimed final estimate.
\end{proof}

\section{The trigonometric polynomial at rationals with small denominator}
\label{SmallAPCM}
In Section~\ref{SmallAPAM}, we analysed $S(\alpha;N)$ when $\alpha$ is close to a 
rational $a/q$ with $q$ being small by using Dirichlet series. The two consequences 
are the lack of effectivity when $q>(\log N)^2$ and the lesser generality of the 
proof. The main result of this section is Lemma~\ref{TrigoSmall}, which restores
effectivity. A closer look discloses that we would be able to accommodate a general sieve setting. We refrain from doing 
so because we want the reader to see clearly where we fail, in this specific example,
to get $\sqrt{q}/\varphi(q)$ as the main term of our bound.

As in Lemma~\ref{CuiWu}, we shall understand the analytic factor
$e(n\alpha)$ in $S(\frac aq+\beta;N)$ as a control of the behaviour in
small intervals. Hence, to handle it, it is enough to understand the
behaviour of $b(n)e(na/q)$ in small intervals. This is the topic of this section. We start with a sieve lemma.
\begin{lem}
\label{auxsieve}
  Assume $2\le z^2\le L$, $m\ge1$ and $\log(2m)\le \sqrt{\log z}$. Then we have 
  \begin{equation*}
      \sum_{\substack{m=m_1-m_2\\ (m_1m_2,P_{4;3}(z))=1\\ M<m_1\le M+L}}1
      \ll \frac{m}{\varphi(m)}\frac{L}{\log z}.
  \end{equation*}
\end{lem}
\begin{proof}
  This is, for instance, a simple instance of Montgomery's sieve. Let us
  denote by $\Sigma$ the quantity to be bounded above. We give some
  details. Notice that $(m_1,P_{4;3}(z))=1$. The mentioned sieve
  provides us with the bound $ \Sigma \le (L+z^2)/{G(z)}$ where
  $G(z)=\sum_{\substack{d\le z\\
      d|P_{4;3}(z))=1}}\mu^2(d)\prod_{p|d}g_m(p)$ and
  \begin{equation*}
    g_m(p)
    =\begin{cases}
        \frac{2}{p-2}&\text{if $(m,p)=1$ and $p|P_{4;3}(z)$,}\\
        \frac{1}{p-1}&\text{if $p|m$ and $p|P_{4;3}(z)$},\\
        0&\text{when $p\nmid P_{4;3}(z)$.}
      \end{cases}
  \end{equation*}
  We appeal for instance to the book \cite[Chapter 2, Theorem 2]{Greaves*01} of G.~Greaves to evaluate the function~$G$. Notice that
  \begin{equation*}
      \sum_{p\le t}g_m(p)\log p
      =\log t+\Ocal(\log\log(3m)).
  \end{equation*}
  We thus infer that
  \begin{equation*}
      G(z)\gg \prod_{\substack{p\le z\\ p|P_{4,;3}(z)}}\bigl(1+g_m(p)\bigr)
      \gg \prod_{\substack{p|m\\ p|P_{4;3}(z)}}\frac{p-2}{p-1}(\log z).
  \end{equation*}
  We finally notice that
  \begin{equation*}
      \prod_{\substack{p|m\\ p|P_{4;3}(z)}}\frac{p-1}{p-2}
      \le \prod_{\substack{p|m\\ p|P_{4;3}(z)}}\frac{p-1}{p-2}
      \frac{p-1}{p}
      \prod_{\substack{p|m\\ p|P_{4;3}(z)}}\frac{p}{p-1}
      \ll m/\varphi(m).
  \end{equation*}
  The lemma follows readily.
\end{proof}

\begin{lem}
  \label{TrigoSmall0}
  Let $A, B\ge1$.
  Let $ N/(\log N)^B\le U\le N$,  $q\le (\log N)^B$ and $(a,q)=1$. For any $A\ge1$, we have
  \begin{equation*}
    \frac{S(\frac{a}{q};N)-S(\frac{a}{q};N-U)}{U/\sqrt{\log N}}
    \ll_{A,B}
    \frac{\sqrt{q\log\log N}\log\log\log N}{\varphi(q)}
    +\frac{1}{(\log N)^A}.
  \end{equation*}
  The implied constant can be effectively computed.
\end{lem}
This can be compared with \cite[Theorem 2a, Chapter IX]{Vinogradov*54}
in the book of I.~M.~Vinogradov, where, in the case of primes, the
author gets $(\log N)^\varepsilon$ (for any positive $\varepsilon$)
rather than our $\sqrt{\log\log N}\log\log\log N$.
\begin{proof}[Proof of Lemma~\ref{TrigoSmall0}]
  We may assume that $N$ is so large that
  $N^{1/2}q\sqrt{\log N}\le N^{3/4}\le U$, for the bound we claim is otherwise trivial to obtain.
  We use Theorem~\ref{Vinobis} with $D=N^{1/4}$. 
\proofsubsection{Estimation of the linear part}
  As in the proof of \eqref{est1}, we readily get
\begin{equation}
  \label{est1s}
  \sum_{\substack{d|P_{4;3}(z)\\
      d\le D\\ q\nmid 4d}}\mu(d)
  \sum_{\substack{N-U<n\le N\\ n\equiv0[d]\\ n\equiv1[4]}}e(na/q)
  \ll \sum_{1\le b\le q}\frac{D}{q}\frac{q}{b}\ll N^{1/4}\log 4q.
\end{equation}
When $q|4d$, we need to refine~\eqref{est1s}. We write $q=q_2q_*$ with $q_2|4$ and $(q_*,2)=1$ and proceed as follows. A reasoning similar to the one used in Lemma~\ref{helper} leads to
\begin{equation}
  \label{est2s}
  \sum_{\substack{d|P_{4;3}(z)\\
      d\le D\\ \; \;  q\mid 4d }}\mu(d)\sum_{\substack{n\equiv0[d]\\ n\equiv
      1[4]\\ N-U<n\le N}}e(na/q)
  =\sum_{\substack{d|P_{4;3}(z)\\
      d\le D\\ q_*| d}}\frac{e(aq_*/q_2)\mu(d)U}{4d}+\Ocal(D/q).
\end{equation}
On using
Lemma~\ref{RTBrun}, we complete the sum and find that
\begin{align}
  \label{est2sb}
  \sum_{\substack{d|P_{4;3}(z)\\
      d\le D\\ q_*| d}}\frac{\mu(d)}{d}
  &=\frac{\mu(q_*)}{q_*}\prod_{p<z, p\equiv
    3[4]}\left(1-\frac{1}{p}\right)\prod_{p\mid
    q_*}\left(1-\frac{1}{p}\right)^{-1}
    \nonumber
  \\&\qquad\qquad+\Ocal\biggl({\sqrt{\log z}\exp\biggl(\frac{-\log D}{\log z}\log\frac{\log D}{2\log z}\biggr)}\biggr)\nonumber\\
 & \ll \frac{1}{\varphi(q)\sqrt{\log z}}+\sqrt{\log z}\exp\biggl(\frac{-\log D}{\log z}\log\frac{\log D}{2\log z}\biggr),
  \end{align}
uniformly in $q\le N$ and provided that $\log z\ge \sqrt{\log N}$ say.

\proofsubsection{Bilinear part -- Preliminary treatment}
Let us now turn towards the treatment of the bilinear part.
We localize the variable $p$ in a diadic interval
$[P,P']\subset [z,\sqrt{N}]$ with $P'\le 2P$.
It turns out that this would be enough when $U=N$ but does not
suffice in general. So we subdivide the diadic interval $[P,P']$ in
subintervals of size at most  $W$: say $[P^*,P^*+W]$, meaning that we
study
\begin{equation*}
  \Sigma_P(W,P^*)
  =\sum_{\substack{N-U<mp\le N\\ (m,P_{4;3}(z))=1\\ P^*\le p\le
      P^*+W\\ p\equiv 3[4]}}
  \mkern-18mu
  \rho(m)e(mpa/q)\;. 
\end{equation*}
In this sum, we replace the condition $mp\in(N-U,N]$ by
$mP^*\in (N-U,N]$, getting a new sum we denote by
$\Sigma^0_P(W,P^*)$. For each $p$, we find that
\begin{equation*}
  \frac{N-U}{P^*}-\frac{N-U}{p}+\frac{N}{P^*}-\frac{N}{p}
  \le \frac{(N-U) W}{P^{*2}}+\frac{N W}{P^{*2}}\le \frac{2N W}{P^2}.
\end{equation*}
Since we have about $\frac12P/\log P$ primes $p$, and as this upper
bound is $\ge1$, this amounts to a
total number of couples $(p,m)$ that is at most $\Ocal(NW/[P\log P])$.
We directly select $W=\sqrt{P}$ and use
$\Sigma^0_P(\sqrt{P}, P^*)=\Sigma^0_P(P^*)$.

\proofsubsection{Bilinear part -- Applying Cauchy}
Let us resume the study of $\Sigma$ by
employing Cauchy's inequality but we keep part of the information that
$\ell$ is prime by using an enveloping sieve up to $P^{\eta}$. This means that we majorize the characteristic function of the primes between $P^*$ and $P^*+W$ by an upper bound sieve. We do so by using the Selberg sieve and add a little twist: we do not sieve out the prime factors of~$q$. The readers may find the details for instance in the paper \cite{van-Lint-Richert*65} by van J. E.~Lint \& H.-E.~Richert.
Here is what we get:
\begin{multline*}
  |\Sigma^0_P(P^*)|^2
  \ll \frac{W}{\log P}\sum_{\substack{ P^*\le \ell \le P^*+W}}\biggr(\sum_{\substack{d|\ell\\ (d,q)=1}}\lambda_d\biggr)^2
    \biggl|
  \sum_{\substack{N-U<mP^*\le N\\ (m,P_{4;3}(z))=1}}
    \rho(m)e(m\ell a/q)
  \biggr|^2
  \\
  \ll \frac{W}{\log P}
  \sum_{\substack{\frac{N-U}{P^*}<m_1,m_2\le \frac{N}{P^*}\\ (m_1m_2,P_{4;3}(z))=1}}
  \mkern-18mu\rho(m_1)\rho(m_2)
  \sum_{\ell\in I(m_1,m_2)}
  \biggl(\sum_{\substack{d|\ell\\ (d,q)=1}}\lambda_d\biggr)^2
  e\Bigl(\frac{(m_1-m_2)\ell a}{q}\Bigr)
\end{multline*}
where $I(m_1,m_2)$ is an interval of length at most $W$. Notice that
$m_1-m_2=m$ has at most $\frac{m}{\varphi(m)}\frac{U}{P\log z}$ solutions when $m\neq0$ by Lemma~\ref{auxsieve}. 
\proofsubsection{Bilinear part -- Diagonal contribution}
Since $m$ lies in an interval of length $\Ocal(U/P)$, the diagonal contribution $C$, i.e. when $m\equiv 0[q]$, satisfies
\begin{equation*}
   C\ll \frac{W}{\log P}
    \biggl({\small\color{red}\text{\footnotesize [m=0]}}\frac{U}{P\sqrt{\log z}}
    +
    \sum_{q|m\ll U/P}
    \frac{m}{\varphi(m)}\frac{U}{P\log z}
    \biggr)
    \frac{W}{\frac{\varphi(q)}{q}\log P}.
\end{equation*}
We write $m=qt$ and notice that $\frac{m}{\varphi(m)}\le \frac{q}{\varphi(q)}\frac{t}{\varphi(t)}$ so that 
\begin{align*}
    C
    &\ll
    \biggl(\frac{U}{P\sqrt{\log z}}
    +
    \frac{U^2}{\varphi(q)P^2\log z}
    \biggr)
    \frac{qW^2}{\varphi(q)(\log P)^2}
    \\&\ll
    \biggl(\frac{qP\sqrt{\log z}}{U}
    +
    1
    \biggr)
    \frac{qU^2W^2}{(\varphi(q)P\log P)^2\log z}
    \ll 
    \frac{qU^2W^2}{(\varphi(q)P\log P)^2\log z}.
\end{align*}
since $P\le N^{1/2}$ and $N^{1/2}q\sqrt{\log N}\le U$.
\proofsubsection{Bilinear part -- Off-diagonal contribution}
By off-diagonal, we mean $m\not\equiv 0[q]$.
First notice that
\begin{align*}
  \sum_{\ell\in I(m_1,m_2)}
  \biggl(\sum_{\substack{d|\ell\\ (d,q)=1}}\lambda_d\biggr)^2
  e\Bigl(\frac{m\ell a}{q}\Bigr)
  &=
    \sum_{\substack{d_1,d_2\le P^{2\eta}\\ (d_1d_2,q)=1}}
  \lambda_{d_1}\lambda_{d_2}
  \sum_{\substack{\ell\in I(m_1,m_2)\\ [d_1,d_2]|\ell}}e\Bigl(\frac{m\ell a}{q}\Bigr)
  \\&\ll qP^{4\eta}.
\end{align*}
The contribution of the $m$-sum is
\begin{equation*}
    \ll \sum_{m\ll U/P}
    \frac{m}{\varphi(m)}\frac{U}{P\log z}
    \ll \frac{U^2}{P^2\log z}.
\end{equation*}
\proofsubsection{Bilinear part -- Resuming the proof}
We gather both estimates and reach
\begin{align*}
  |\Sigma^0_P(P^*)|^2
  &\ll
  \frac{qU^2W^2}{(\varphi(q)P\log P)^2\log z}
  +
  qP^{4\eta}\frac{WU^2}{P^2(\log P)\log z}
  \\&\ll
   \frac{qU^2W^2}{(\varphi(q)P\log P)^2\log z}
 \biggl(
  1
  +
  \frac{\varphi(q)^2\log P}{W}
  P^{4\eta}
  \biggr).
\end{align*}
We select $\eta=1/9$ and assume that $q\le z^{1/8}$, so that
\begin{equation*}
    |\Sigma^0_P(P^*)|^2
    \ll 
    \frac{qU^2W^2}{(\varphi(q)P\log P)^2\log z}
    \qquad(q\le z^{1/8}).
\end{equation*}
We sum over $P^*$, recall that $W=\sqrt{P}$. This leads to
\begin{equation*}
  \sum_{P^*}|\Sigma^0_P(P^*)|
  \ll
  \frac{\sqrt{q}U}{\varphi(q)(\log z)^{1/2}\log P}
  .
\end{equation*}
We add the error term $\Ocal(NP^{1/2}/[P\log P])$
and then sum over $P$, add~\eqref{est1s},~\eqref{est2s} and~\eqref{est2sb}
together with the error term coming from Theorem~\ref{Vinobis},
and, on calling $S^\sharp(a/q)$ the trigonometric sum to be
studied, we find that
\begin{multline*}
  S^\sharp(a/q)
  \ll
  \frac{\sqrt{q}U\log\frac{\log N}{\log z}}{(\log z)^{1/2}\varphi(q)}
  + {\color{red}\text{\footnotesize [$P$ to $P^*$]}}
  \frac{\sqrt{N}}{\sqrt{z\log z}}
  \\+
  {\color{red}\text{\footnotesize [\eqref{est2s}, \eqref{est2sb}]}}\frac{U}{\varphi(q)\sqrt{\log z}}
  +{\color{red}\text{\footnotesize [\eqref{est1s}]}}N^{1/4}\log 4q
  \\+ {\color{red}\text{\footnotesize [$E'$ from Thm~\ref{Vinobis}]}}
  \frac{U}{z}
  +N\exp\biggl(-\frac{\log N}{4\log z}\log\frac{\log N}{8\log
    z}\biggr)
  \sqrt{\log z}.
\end{multline*}
Some explanation may be required for the second summand of the first line.
In order to sum $1/\log P$ over our diadic intervals delimited by
$P=z2^k$, where
$k\in\{0,\ldots, K\}$ with $K=[(\log \frac{\sqrt{N}}{z})/\log 2]$, we
first notice that (with $T=[(\log z)/\log 2]$)
\begin{align*}
  \sum_{0\le k\le K}\frac{1}{\log z+k\log 2}
  &\le \frac{1}{\log 2}\sum_{0\le k\le K}\frac{1}{k + T}
  \\&\ll \sum_{ k\le K+T}\frac{1}{k}-\sum_{k<T}\frac{1}{k}
  \ll \log\frac{\log N}{\log z}+\frac{1}{\log z}
  \\&\ll \log\frac{\log N}{\log z}.
\end{align*}
Since $q\le \sqrt{N}$, we easily simplify the above in
\begin{equation*}
  S^\sharp(\alpha)
  \ll
  \frac{U\sqrt{q}\log\frac{\log N}{\log z}}{\varphi(q)\sqrt{\log z}}
  +N^{1/2}
  +\frac{U}{z}
  +N\exp\biggl(-\frac{\log N}{4\log z}\log\frac{\log N}{8\log
    z}\biggr)
  \sqrt{\log z}.
\end{equation*}
We select $\log z=\frac{\log N}{(4A+5)\log\log N}$ and reach
\begin{align*}
   S^\sharp(\alpha)
  &\ll
  \frac{U\sqrt{q\log\log N}\log\log\log N}{\varphi(q)\sqrt{\log N}}
  +N\exp\bigl(-(A+1)(\log\log N)\log\log\log N\bigr)
  \\&\ll
  \frac{U\sqrt{q\log\log N}\log\log\log N}{\varphi(q)\sqrt{\log N}}
  +\frac{N}{(\log N)^{(A+1)\log\log\log N}}.
\end{align*}
This estimate would allow $q\ll (\log N)^{C\log\log\log N}$ for a
positive constant~$C$ and similarly $U\ge N (\log N)^{-C\log\log\log
  N}$. We prefer to simplify the statement. The proof of this lemma is ended.
\end{proof}

Let us convert the behaviour in small intervals in analytic variations
\begin{lem}
  \label{TrigoSmall}
  Let $A,B\ge1$.
  Let $q\le (\log N)^B$, $(a,q)=1$ and $|\alpha-\frac aq|\le(\log N)^B/N$. We have
  \begin{equation*}
    \frac{S(\alpha;N)}{N/\sqrt{\log N}}
    \ll_{A,B}
    \frac{\sqrt{q\log\log N}\log\log\log N}{\varphi(q)}
    +\frac{1}{(\log N)^A}.
  \end{equation*}
  The implied constant can be effectively computed.
\end{lem}

\begin{proof}
  Let us write $\alpha=\frac aq+\beta$. We write
  \begin{equation*}
    S(\alpha;N)
    =\sum_{1\le k\le \frac{N}{U}}\bigl(S(\alpha;N-(k-1)U)-S(\alpha;N-kU)\bigr)
    +\Ocal\biggl(\frac{U}{\sqrt{\log U}}\biggr)
    \end{equation*}
    for a parameter $U$ that we shall choose later. Given $k$ and on
    setting $M=N-(k-1)U$, we have
\begin{align*}
        S(\alpha;M)-&S(\alpha;M-U)
      \\=
        &e(M\beta)\sum_{M-U<n\le M}b(n)e(na/q)
      -2\pi i\beta 
      \int_{M-U}^{M}\sum_{M-U<n\le t}b(n)e(na/q)e(t\beta)dt 
      \\&
      \ll (1+|\beta\; U|)\biggl(\frac{U\sqrt{q}(\log\log M)^{1/2}\log\log\log M}{(\log M)^{1/2}\varphi(q)}\biggr)
    \end{align*}
    by Lemma~\ref{TrigoSmall0}. 
    On summing over $k$, we get
    \begin{equation*}
      S(\alpha;N)
      \ll (1+|\beta|U)\biggl(
      \frac{\sqrt{q}N(\log \log N)^{1/2}\log\log\log N}{\varphi(q)\sqrt{\log N}}
      +
      \frac{N}{(\log N)^A}\biggr).
    \end{equation*}
    We select $U=N/(\log N)^B$ and conclude easily.
\end{proof}

\section{A global bound for the trigonometric polynomial, proof of Theorem~\ref{Trigo}}
\label{ProofTrigo}
\begin{proof}[Proof of Theorem~\ref{Trigo}]
  Let us start with the generic non effective part.
  Lemma~\ref{TrigoLargeII} applies when $q$ is larger than~$(\log N)^{14+2A}$,
  while Lemma~\ref{mainexpter} is enough to take care of the case when
  $q\le (\log N)^{14+2A}$.

  Let us now examine how to modify the previous sketch to recover
  effectivity when $A<1/2$. Lemma~\ref{TrigoLargeII} applies when $q\ge(\log N)^{14+2A}$, Lemma~\ref{TrigoSmall} applies when
  $\log N\le q\le (\log N)^{14+2A}$ and Lemma~\ref{mainexpter} is enough
  to take care of the case $q\le \log N$.
\end{proof}

\section{An additive application, Proof of Theorem~\ref{th2}}
\label{Proofth2}
\begin{lem}\label{L1}
  For any $\ell>2$ and any subset $\mathcal{B}'$ of the Gaussian integers, we have
  \begin{equation*}
    \int_0^1\biggl|
    \sum_{\substack{b\le N\\ b\in\mathcal{B} '}}e(b\alpha)\biggr|^\ell
    d\alpha\ll \frac{N^{\ell-1}}{(\log N)^{\ell/2}}.
  \end{equation*}
\end{lem}
\begin{proof}
This is a trivial consequence of the method developed
in~\cite[Theorem 1.1]{Green-Tao*04} by B.~Green and T.~Tao, or in a
simplified form in~\cite{Ramare*22}.
\end{proof}

Here is the main lemma of this section.
\begin{lem}
  \label{inith2}
  Let $N\equiv3[4]$ be an integer, $K\in[2, \log\log N]$ be a
  parameter, and $\mathcal{B}_1$ and $\mathcal{B}_2$ be two sets of
  odd primitive Gaussian integers. We have
  \begin{multline*}
    \sum_{\substack{n+b_1+b_2=N\\ b_1\in\mathcal{B}_1, b_2\in\mathcal{B}_2}}\mkern-10mub(n)
    =
    \frac{4CP_{4;3}(K)}{\varphi(P_{4;3}(K))}
    \mkern-20mu\sum_{\substack{b_1\in\mathcal{B}_1,
        b_2\in\mathcal{B}_2\\ N-(b_1+b_2)\ge 3\\
        (N-(b_1+b_2),P_{4,3}(K))=1}}
    \mkern-10mu
    \frac{1}{\sqrt{\log(N-(b_1+b_2))}}
    \\[-10pt]+\Ocal\biggl(
    \frac{N^2}{\sqrt{K}(\log N)^{3/2}}
    \biggr).
  \end{multline*}
\end{lem}

\begin{proof}
We consider the Dirichlet dissection of the torus
$\mathbb{R}/\mathbb{Z}$:
\begin{equation*}
  \mathfrak{M}(a/q)=\{\alpha / |\alpha-a/q|\le 1/(qQ)\}
\end{equation*}
for $Q=N/(\log N)^{B}$, where $B=2A+14$ for some $A>0$.  We call each interval $\mathfrak{M}(a/q)$ as a major arc and we denote their union as 
\begin{equation}
    \label{majorarcs}
    \mathfrak{M}= \bigcup_{q_2\in\{1,2,4\}}
  \bigcup_{\substack{q_*\le
      N/Qq_2\\ q_*| P_{4;3}(K)}}
      \bigcup_{a\mode q_2q_*}\mathfrak{M}\Bigl(\frac{a}{q_2q_*}\Bigr)\; .
\end{equation}

 We set
\begin{equation*}
  \mathfrak{m}=\bigl(\mathbb{R}/\mathbb{Z}\bigr)\setminus
  \mathfrak{M}
      .\end{equation*}
     
We further define 
\begin{equation}
  S_1(\alpha)=\sum_{\substack{b_1\le N\\ b_1\in \mathcal{B}_1}}
  e(b_1\alpha),
  \quad
  S_2(\alpha)=\sum_{\substack{b_2\le N\\ b_1\in\mathcal{B}_2}}
  e(b_2\alpha),
\end{equation}
while $S^\flat_K(\alpha;N)$ is defined in \eqref{defSflatK}.
We employ the circle method:
\begin{equation}
  \label{def:rofn}\left\{
    \begin{aligned}
      r(N)&=\int_0^1 S(\alpha)S_1(\alpha)S_2(\alpha)e(-N\alpha)d\alpha,
      \\r^\flat(N)&=\int_0^1
      S^\flat_K(\alpha,N)S_1(\alpha)S_2(\alpha)e(-N\alpha)d\alpha,
    \end{aligned}
  \right.
\end{equation}
and write
\begin{equation}
  r(N)=r_0(N)+r'(N),\quad
  r'(N)=\int_{\mathfrak{m}}S(\alpha)S_1(\alpha)S_2(\alpha)e(-N\alpha)d\alpha
\end{equation}
and, correspondingly,
\begin{equation}
  r^\flat(N)=r^\flat_0(N)+r^{\flat\prime}(N),\quad
  r^{\flat\prime}(N)=
  \int_{\mathfrak{m}}S^\flat_K(\alpha)S_1(\alpha)S_2(\alpha)e(-N\alpha)d\alpha.
\end{equation}
Let $\ell\in(2,3)$. We notice that
\begin{equation}
  \label{eq:3}
  \frac{1}{\ell}+\frac{1}{\ell}+\frac{\ell-2}{\ell}=1,\quad
  \frac{\ell}{\ell-2}=\ell+\frac{3-\ell}{\ell-2}\ell
\end{equation}
and use H\"older inequality to infer that
\begin{equation*}
  |r'(N)|\le
  \biggl(\int_{\mathfrak{m}}|S_1(\alpha)|^{\ell}d\alpha\biggr)^{1/\ell}
  \biggl(\int_{\mathfrak{m}}|S_2(\alpha)|^{\ell}d\alpha\biggr)^{1/\ell}
  \biggl(\int_{\mathfrak{m}}|S(\alpha)|^{\ell}
  |S( \alpha)|^{\frac{\ell(3-\ell)}{\ell-2}}d\alpha\biggr)^{\frac{\ell-2}{\ell}}.
\end{equation*}
We extract $\max_{\alpha\in\mathfrak{m}} |S(
\alpha)|^{\frac{\ell(3-\ell)}{\ell-2}}$ from the last term and then
extend all three integrals to the full circle. When $\alpha\in \mathfrak{m}$, we have three cases: (i) $q\geq \log^{B}(N)$  (ii) $q\leq  \log^{B}(N), q_2\nmid 4$ (iii)  $q\leq  \log^{B}(N), q_2\mid 4, q_*\nmid P_{4,3}(K)$. We majorize the maximum in case (i)  via Theorem~\ref{Trigo} and via Lemma~\ref{mainexpter} in the other cases where we further notice that $\log q\ll \log\log N$. In the case (iii), there exists at least one prime $p>K$ that divides $q^*$ and thus $1/\varphi(q^*)\ll 1/K$. This argument implies that
\begin{equation*}
 \max_{\alpha\in\mathfrak{m}} |S(\alpha)|
 \ll \max\biggl(\frac{N}{(\log N)^A}, 
 \frac{N\log\log N}{(\log N)^{3/2}},
 \frac{N}{K\sqrt{\log N}}\biggr) \ll \frac{N}{K\sqrt{\log N}}\; .
\end{equation*}
Hence we get that 
\begin{equation*}
 |r'(N)|\ll \biggl(\frac{N^{\ell-1}}{\log^{\ell/2}(N)}\biggr)\biggl(\frac{N}{K\sqrt{\log N}}\biggr)^{3-\ell}\ll \frac{N^2}{K^{3-\ell}(\log N)^{3/2}}\; .
\end{equation*}
We choose $\ell=5/2$.

The same proof applies for $r^{\flat\prime}(N)$, though the estimates for
$S^\flat_K(\alpha)$ are different; furthermore Lemma~\ref{L1b}
replaces Lemma~\ref{L1}. We obtain:
\begin{equation*}
  r^{\flat\prime}(N)\ll 
  \biggl(e^K\frac{N}{(\log N)^B}\biggr)^{3-\ell}
  \biggl(\frac{N^{\ell-1}}{(\log N)^ {\ell/2}}\biggr)^{\frac{\ell-2}{\ell}+\frac{2}{\ell}} \biggl(\sqrt{\log K}\biggr)^{\frac{\ell-2}{\ell}}
\end{equation*}
as $q\le N/(\log N)^B$. This gets reduced to
\begin{equation*}
  r^{\flat\prime}(N)\ll \frac{N^2(\log N)^{3-\ell}\sqrt{\log\log N}}{(\log N)^{B(3-\ell)+\frac{\ell}{2}}}
  \ll \frac{N^2}{(\log N)^{B-1}}.
\end{equation*}
On the major arcs, $S(\alpha;N)$ and $S^\flat_K(\alpha;N)$ are
similar by Lemma~\ref{mainexpter} and~\ref{model}, hence the result. 
\end{proof}

\begin{proof}[Proof of Theorem~\ref{th2}]
  We start with Lemma~\ref{inith2} with the aim of replacing the
  $1/\sqrt{\log(N-b_1-b_2)}$ by the simpler $1/\sqrt{\log N}$.
  Let $\varepsilon=1/(\log N)^2$.
  We first notice that, with $m=b_1+b_2$,
  \begin{equation*}
    \frac{1}{\sqrt{\log(N-m)}}-\frac{1}{\sqrt{\log N}}
    =
    \frac{\log\frac{N}{N-m}}{\sqrt{\log(N-m)\log
        N}(\sqrt{\log N}+\sqrt{\log(N-m)}}.
  \end{equation*}
  Hence
  \begin{multline*}
    \sum_{\substack{b_1\in\mathcal{B}_1,
        b_2\in\mathcal{B}_2\\ N-(b_1+b_2)\ge \varepsilon N\\
        (N-(b_1+b_2),P_{4,3}(K))=1}}
    \frac{1}{\sqrt{\log(N-b_1-b_2)}}
    =
    \sum_{\substack{b_1\in\mathcal{B}_1,
        b_2\in\mathcal{B}_2\\ N-(b_1+b_2)\ge \varepsilon N\\
        (N-(b_1+b_2),P_{4,3}(K))=1}}
    \frac{1}{\sqrt{\log N }}
    \\+\Ocal\biggl(
    \frac{|\mathcal{B}_1\cap[1,N]||\mathcal{B}_2\cap[1,N]|\log(1/\varepsilon)}
    {\sqrt{\log(\varepsilon N)}\log N}\biggr).
  \end{multline*}
  This remainder term is readily seen to be $\ll N^2(\log\log N)/(\log
  N)^{5/2}$. Next, there for each integer $m$ in $[(1-\varepsilon)N,N]$,
  the sieve tells us that number of decomposition $m=b_1+b_2$ is at most $N(\log\log
  N)/\log N$.
  This gives us finally that
  \begin{multline*}
    \sum_{\substack{b_1\in\mathcal{B}_1,
        b_2\in\mathcal{B}_2\\ N-(b_1+b_2)\ge 3\\
        (N-(b_1+b_2),P_{4,3}(K))=1}}
    \mkern-10mu
    \frac{1}{\sqrt{\log(N-(b_1+b_2))}}
    =\sum_{\substack{b_1\in\mathcal{B}_1,
        b_2\in\mathcal{B}_2\\ N-(b_1+b_2)\ge 3\\
        (N-(b_1+b_2),P_{4,3}(K))=1}}
    \mkern-10mu
    \frac{1}{\sqrt{\log N}}
    \\+
    \Ocal\biggl(
    \frac{N^2\log\log N}{(\log N)^{5/2}}
    \biggr).
  \end{multline*}
  Theorem~\ref{th2} follows readily for this estimate.
\end{proof}

\section{A challenge. Proof of Theorem~\ref{ValAtSqrtOf2}}
\label{ProofChallenge}
The following is useful for the proof of Theorem~\ref{ValAtSqrtOf2}.
\begin{lem}
  \label{ConVSqrt2}
  For any integer $k$, we define 
  $ N_k=\frac{(1+\sqrt{2})^k-(1-\sqrt{2})^k}{2\sqrt{2}}$.
  Then the  quantity $N_k$ is an integer, $N_k$ and  $N_{k+1}$ are
  coprime, and we have
  $
    \bigl|\sqrt{2}-\frac{N_{k+1}-N_k}{N_k}\bigr|\le \frac{1}{2N_k^2}.
    $
  When $y$ is large enough, each interval $[y,3y]$ contains a $N_k$.
\end{lem}

\begin{proof}
  We readily check that  $N_1=1$, $N_2=2$ and that
  $N_{k+2}-2N_{k+1}-N_k=0$, from which the reader will easily deduce
  that the $N_k$'s are integers and that $N_k$ and $N_{k+1}$ are
  coprime. Let us set $\alpha=1+\sqrt{2}$ so that
  $1/\alpha=-(1-\sqrt{2})$ and $\alpha^2=2\alpha+1$. We compute that
  \begin{equation*}
    8N_k^2\biggl(\sqrt{2}-\frac{N_{k+1}-N_k}{N_k}\biggr)
    =8\alpha N_k^2-8N_kN_{k+1}
    = 2(\alpha-1)(\alpha^{-2k}-(-1)^k)
  \end{equation*}
  and this last quantity is swiftly shown to be, in absolute value, at
  most~4. The last assertion comes from the fact that $1+\sqrt{2}< 3$.
\end{proof}

\begin{proof}[Proof of Theorem~\ref{ValAtSqrtOf2}]
  This proof is a direct consequence of two points:
  Lemma~\ref{TrigoLargeII} and the fact that $\sqrt{2}$ has convergent
  denominators in each interval $[y,3y]$, as shown in Lemma~\ref{ConVSqrt2}.
\end{proof}
\section{Bounds on multiples. Proof of Theorem~\ref{Family}}
\label{ProofFamily}
\begin{proof}[Proof of Theorem~\ref{Family}]
Let us assume that $N$ is not the square of an integer.
  We start with Theorem~\ref{AS} which we use with $Z=\sqrt{N}$, $u_n=e(r\alpha n)$, $M=M_0$ and the parameters $M_0$ and $z$ to be chosen.
  
  \proofsubsection{Bilinear part -- Vinogradov part}
  We are to study
  \begin{equation*}
      \Delta_1=\sum_{\substack{z\le p\le \sqrt{N}\\ p\equiv3[4]}}
      \sum_{R/2<r\le R}\biggl|
      \sum_{N/2<pm\le N}\rho(m)e(pmr\alpha)\biggr|.
  \end{equation*}
  Let us localize the variable $p$ in $(P',P]$, so the variable $m$ is localized in $(N/(2P),N/P']$. 
  We use Lemma~\ref{HarmanLem4} to infer that this is bounded above in absolute value, and up to a multiplicative constant, by
  \begin{equation*}
      \frac{NR}{\sqrt{\log N}}
      \biggl(
      \frac{1}{\sqrt{q}}
      +\frac{\sqrt{q}}{\sqrt{RN}}
      +
        \frac{\sqrt{P}}{\sqrt{N}}
        + \frac{1}{\sqrt{RP}}
      \biggr)
      (\log N)^{4}.
  \end{equation*}
  On summing diadically over $P'$, we infer that
  \begin{equation}
  \label{BoundDelta1}
      \Delta_1\ll 
      \frac{NR}{\sqrt{\log N}}
      \biggl(
      \frac{1}{\sqrt{q}}
      +\frac{\sqrt{q}}{\sqrt{RN}}
      +
        \frac{1}{N^{1/4}}
        + \frac{1}{\sqrt{Rz}}
      \biggr)
      (\log N)^{5}.
  \end{equation}
  \proofsubsection{Linear part}
  We are to study
  \begin{equation*}
      \Delta_0=
      \sum_{R/2<r\le R}c_r
      \sum_{\substack{d\le M_0\\ d|P_{4;3}(z)}}\mu(d)
      \sum_{\substack{n\le N\\ n\equiv 0[d]}}e(nr\alpha).
  \end{equation*}
  Let us glue $r$ and $d$ in a single variable $u$ and use Lemma~\ref{VaughanLem3} with $U=RM_0$ and $UV=N$, getting, for any positive $\ve$,
  \begin{equation}
     \label{BoundDelta0}
      \Delta_0\ll_\ve
      \biggl(
      \frac{N}{q}+RM_0+q
      \biggr)
      N^\ve.
  \end{equation}
  
  \proofsubsection{Bilinear part -- Simple sieve part}
  We are finally to study, for any $t\in[0,1]$ the quantity
  \begin{equation*}
      \Delta_3
      =
      \sum_{R/2<r\le R}\biggl|
      \sum_{M_0\le \ell\le M_0z}
      \alpha_\ell(t)
      \sum_{N/2 < k\ell \le N}\beta_k(t)e(k\ell r\alpha)
      \biggr|.
  \end{equation*}
  Let us localize the variable $\ell$ in an interval $(L',L]$, so that $k$ is localized on $(N/(2L), N/L']$. We may thus employ Lemma~\ref{HarmanLem4} and reach the bound, for each of any such piece and up to a multiplicative constant,
  \begin{equation*}
      \frac{RN}{\sqrt{\log N}}
      \biggl(
      \frac{1}{\sqrt{q}}
      +\frac{\sqrt{q}}{\sqrt{RN}}
      +\frac{1}{\sqrt{N/L}}
      +\frac{1}{\sqrt{RL}}
      \biggr)
      (\log N)^5.
  \end{equation*}
  We sum diadically on $L$ to reach
  \begin{equation}
      \label{BoundDelta3}
      \Delta_3\ll
      \frac{RN}{\sqrt{\log N}}
      \biggl(
      \frac{1}{\sqrt{q}}
      +\frac{\sqrt{q}}{\sqrt{RN}}
      +\frac{\sqrt{M_0z}}{\sqrt{N}}
      +\frac{1}{\sqrt{RM_0}}
      \biggr)
      (\log N)^6.
  \end{equation}
  \proofsubsection{Gathering the estimates}
   It is best to choose $M_0$ as small as possible though larger than $z$, so we select $M_0=z$, getting the upper bound
   \begin{equation*}
       \frac{NR}{\sqrt{\log N}}
      \biggl(
      \frac{1}{\sqrt{q}}
      +\frac{\sqrt{q}}{\sqrt{RN}}
      +
        \frac{1}{N^{1/4}}
        + \frac{1}{\sqrt{Rz}}
        +\frac{z}{\sqrt{N}}
        +\frac{1}{Rq}+\frac{z}{N}+\frac{q}{RN}
        +\frac{1}{z}
      \biggr)
      N^\ve.
   \end{equation*}
   We readily simplify that into
   \begin{equation*}
       \frac{NR}{\sqrt{\log N}}
      \biggl(
      \frac{1}{\sqrt{q}}
      +\frac{\sqrt{q}}{\sqrt{RN}}
      +
        \frac{1}{N^{1/4}}
        + \frac{1}{\sqrt{Rz}}
        +\frac{z}{\sqrt{N}}
        +\frac{1}{z}
      \biggr)
      N^\ve.
   \end{equation*}
   We assume that $R\le N$ and select $z=(N/R)^{1/3}$. Hence the above quantity is bounded by 
   \begin{equation*}
       \frac{NR}{\sqrt{\log N}}
      \biggl(
      \frac{1}{\sqrt{q}}
      +\frac{\sqrt{q}}{\sqrt{RN}}
      +
        \frac{1}{N^{1/4}}
        +\frac{1}{R^{1/3}N^{1/6}}
        +\frac{R^{1/3}}{N^{1/3}}
      \biggr)
      N^\ve.
   \end{equation*}
   We readily check that ${R^{-1/3}N^{-1/6}}+{R^{1/3}}{N^{-1/3}}\ge N^{-1/4}$, so we may discard this summand.
   This concludes the proof of our lemma.
\end{proof}

\section{Diophantine approximation and proof of Theorem~\ref{Irrat}}
\label{ProofIrrat}
We recall a classical lemma of trigonometric approximation in the precise form given in \cite[Lemma~1]{Harman*83} by G.~Harman which is based on \cite{Montgomery*78}, see also \cite[Theorem 19]{Vaaler*85} by J.~Vaaler.
\begin{lem}
  \label{TrigoApprox} Let $(a_n)$ be a sequence of non-negative reals for $N/2<n\le N$. Then for $\delta>0$, $\alpha$ and $\beta$ be arbitrary real numbers and $R\ge1$ be some given parameter, we have
  \begin{equation*}
      \sum_{\|n\alpha-\beta\|\le \delta}a_n
      =2\delta \sum_n a_n
      +\Ocal\biggl(\delta \sum_{r\le R}\biggl|\sum_n a_ne(rn\alpha)\biggr|\biggr)
      +\Ocal\biggl(\frac{1}{R}\sum_{n}a_n\biggr).
  \end{equation*}
\end{lem}

\begin{proof}[Proof of Theorem~\ref{Irrat}]
  We employ Lemma~\ref{TrigoApprox} together with Theorem~\ref{Family}. With an obvious notation, this leads to
  \begin{align*}
  \sum_{\|n\alpha-\beta\|\le \delta}b(n)
      &=(2\delta+\Ocal(1/R)) B(N)
      +
  \delta\sum_{r\le R}|S(r\alpha;N)|
  \\&= (2\delta+\Ocal(1/R)) B(N)
      \\&\qquad+\Ocal_\ve\biggl(
      \frac{\delta RN}{\sqrt{\log N}}
      \biggl(
      \frac{1}{\sqrt{q}}
      +\frac{\sqrt{q}}{\sqrt{RN}}
        +\frac{1}{R^{1/3}N^{1/6}}
        +\frac{R^{1/3}}{N^{1/3}}
      \biggr)
      N^\ve\biggr).
      \end{align*}
      Since $\alpha$ is quadratic irrational, there exists a constant $C$ such that we may find approximations $a/q$ of $\alpha$ with $q\in[\sqrt{RN},C\sqrt{RN}]$ and $|q\alpha-a|\le 1/q$ for any $N$, see Lemma~\ref{ConVSqrt2} for the case $\alpha=\sqrt{2}$. With such a choice, we get
      \begin{equation*}
      \biggl|
      \sum_{\|n\alpha-\beta\|\le \delta}\frac{b(n)}{B(N)}-2\delta
      \biggr|
      \ll_\ve\frac{1}{R}
      +
      \delta
      \biggl(
      \frac{\sqrt{R}}{\sqrt{N}}
      +
      \frac{R^{3/4}}{N^{1/4}}
        +\frac{R^{2/3}}{N^{1/6}}
        +\frac{R^{4/3}}{N^{1/3}}
      \biggr)
      N^\ve.
      \end{equation*}
      We select $R=N^{\frac14-2\theta}$ and $\delta^{-1}=N^{\frac{1}{4}-\theta}$ for $\theta>0$ and small enough. This proves our theorem.
\end{proof}

\section{Appendix: some estimates with primes congruent to 3 modulo~4}
We shall freely use the Mertens Theorems for primes congruent to~3
modulo~4, namely
\begin{equation}
  \label{eq:8}
  \sum_{\substack{n\le N\\ n\equiv 3[4]}}\frac{\Lambda(n)}{n}
  =\tfrac12\log N+ \Ocal(1)
  \quad\text{and}\quad
  \prod_{\substack{p\le z\\ p\equiv 3[4]}}\biggl(1-\frac1p\biggr)
  \asymp 1/\sqrt{\log z}.
\end{equation}
An explicit version of the first estimate may be found in
\cite[Corollaire 1]{Ramare*02}, but we wave this aspect most of the times.
However, in the case of Theorem~\ref{Vinobis}, it is
better to rely on the next explicit estimate that is rather
straightforward to get.
\begin{lem}
  \label{mertensp34}
  When $z\ge1$, we have
  $\displaystyle
    \sum_{\substack{p<z\\ p\equiv 3[4]}}\frac{\log p}{p+1}\le\tfrac12\log z-\tfrac14\log 3$.
\end{lem}

\begin{proof}
  The inequality has been checked by GP/Pari up to $10^8$ with a
  minimal value reached for $z=3$.
  Let us proceed for the proof. We first check that
  \begin{equation*}
    \sum_{\substack{p<z\\ p\equiv 3[4]}}\frac{\log p}{p+1}
    \le
    \sum_{\substack{n<z\\ n\equiv 3[4]}}\frac{\Lambda(n)}{n}
    .
  \end{equation*}
  We treat this last sum with Dirichlet characters
  and write
  \begin{equation*}
    2\sum_{\substack{n<z\\ n\equiv 3[4]}}\frac{\Lambda(n)}{n}
    =\sum_{\substack{n<z\\ n\equiv 1[2]}}\frac{\Lambda(n)}{n}
    -\sum_{\substack{n<z}}\frac{\Lambda(n)\chi_4(n)}{n}
  \end{equation*}
  where $\chi_4$ is the only non-principal character modulo~4.
  By using the Prime Number Theorem in Arithmetic Progression, we
  conclude that
  \begin{equation*}
    2\sum_{\substack{n<z\\ n\equiv 3[4]}}\frac{\Lambda(n)}{n}
    =\log z-\log 2 -\gamma +\frac{L'(1,\chi_4)}{L(1,\chi_4)}+o(1),
  \end{equation*}
  see also \cite[Corollaire 1]{Ramare*02}. We have $L(1,\chi_4)=\pi/4$
  (see for instance \cite[Th\'eor\`eme~3, Part~2, Chapter~V]{Borevitch-Chafarevitch*93} by Z.~I.~Borevitch and I.~R.~Chafarevitch) and $L'(1,\chi_4)$ may be computed approximately with the series
  expansion:
  \begin{align*}
    -L'(1,\chi_4)
    &=
      \sum_{k\ge0}\biggl(\frac{\log(4k+1)}{4k+1}-\frac{\log(4k+3)}{4k+3}\biggr)
    \\&=
    \sum_{k\ge0}\frac{2\log(4k+1)-(4k+1)\log(1+\frac2{4k+1})}{(4k+1)(4k+3)}
    \\&=
    -\frac{\log 3}{3}+\sum_{k\ge1}\frac{2\log(4k+1)-(4k+1)\log(1+\frac2{4k+1})}{(4k+1)(4k+3)}.
  \end{align*}
  The good part of the last series is that it is a sum of non-negative
  terms and can thus be computed with acceleration of convergence
  process. We get in this manner
  \begin{equation}
    \label{eq:9}
    L'(1,\chi_4)=0.192901\cdots,\quad
    \frac12\biggl(\frac{L'(1,\chi_4)}{L(1,\chi_4)}-\gamma-\log 2\biggr)
    =-0.512376\cdots\; .
  \end{equation}
  By \cite[Corollaire 4]{Ramare*02}, we have, when $z\ge 182$
  \begin{equation}
    \label{eq:10}
    \sum_{\substack{n<z\\ n\equiv 3[4]}}\frac{\Lambda(n)}{n}
    =\tfrac12\log
    z+\frac12\biggl(\frac{L'(1,\chi_4)}{L(1,\chi_4)}-\gamma-\log
    2\biggr)
    +\Ocal^*(2/9).
  \end{equation}
  On collecting our estimates, the lemma follows readily.
\end{proof}

\printbibliography

\end{document}